%% file: MacPherson.tex
\documentclass{article}
\usepackage{amssymb,amscd,amsmath, amsthm}
\usepackage{epsfig}

\def\star{\mathop{\rm star}}

\def\link{\mathop{\rm link}}
\def\rank{\mathop{\rm rank}}
\def\sign{\mathop{\rm sign}}

\def\MacP{\mathop{\rm MacP}}

\def\Im{\mathop{\rm Im}}

\def\VM{{\cal V}^*(M)}         
\def\Z2{{\mathbb Z}_2}              
\def\Rn{\mathbb R^n}
\def\Rk{\mathbb R^k}
\def\Rinfty{\mathbb R^\infty}

\def\RPinfty{\mathbb R P^\infty}
\def\R{\mathbb R}

\def\Aa{{\cal A}}

\def\Ff{{\cal F}}
\def\Gg{{\Gamma}}

\def\Ll{{\cal L}}
\def\Mm{{\cal M}}
\def\Oo{{\cal O}}
\def\Qq{{\cal Q}}

\def\Vv{{\cal V}}

\def\Zz{\mathbb Z}

\newtheorem{theorem}{Theorem}[section]
\newtheorem{lemma}[theorem]{Lemma}
\newtheorem{prop}[theorem]{Proposition}

\newtheorem{corollary}[theorem]{Corollary}
\newtheorem*{A}{Theorem A}
\newtheorem*{B}{Theorem B}
\newtheorem*{C}{Theorem C}
\newtheorem*{D}{Theorem D}

\newtheorem*{TRP}{Topological Representation Theorem}
\newtheorem*{OHL}{Order Homotopy Lemma}
\newtheorem*{QTA}{Quillen's Theorem A}
\newtheorem*{QTB}{Quillen's Theorem B}
\newtheorem*{BC}{Babson's Criterion}
\newtheorem*{SQF}{Spherical Quasifibration Theorem}
\newtheorem*{CT}{Comparison Theorem}
\newtheorem*{TI}{Thom Isomorphism Theorem}
\theoremstyle{definition}
\newtheorem{defn}[theorem]{Definition}

\newtheorem{remark}[theorem]{Remark}
\newtheorem{example}[theorem]{Example}

\def\beginwval#1#2{\bgroup
  \edef\@savecount{\the\value{#1}}
  \expandafter\def\csname the#1\endcsname{#2}
  \begin{#1}}
\def\endwval#1{\setcounter{#1}{\@savecount}\end{#1}\egroup}

\begin{document}
\title{\mbox{Mod 2 Cohomology of Combinatorial Grassmannians }}
\author{Laura Anderson\footnote{Partially supported by grants from the
National Science Foundation.} \and James F. Davis$^*$}
\date{}
\maketitle

Oriented matroids  have long been of use in various areas of combinatorics \cite{BLSWZ}.
Gelfand and MacPherson \cite{GM} initiated the use of oriented matroids in manifold and
bundle theory, using them to formulate a combinatorial formula for the  rational Pontrjagin
classes of a differentiable manifold.  MacPherson \cite{Mac} abstracted this into a manifold
theory ({\em combinatorial differential (CD) manifolds}) and a bundle theory (which we call {\em
combinatorial vector bundles} or {\em matroid bundles}). In this paper
we explore the relationship between combinatorial vector bundles and real vector bundles.
As a consequence of our results we get theorems relating the topology of the  \emph{combinatorial
Grassmannians} to that of their real analogs.

The theory of oriented matroids gives a combinatorial abstraction of linear algebra; a
$k$-dimensional subspace of $\Rn$ determines a rank $k$ oriented matroid with elements
$\{1,2, \dots, n\}$.  Such oriented matroids can be given a partial order by using the
notion of \emph{weak maps}, which geometrically corresponds to moving the $k$-plane into more special
position with respect to the standard basis of $\Rn$.  The poset $\MacP(k,n)$ of rank $k$
oriented matroids with $n$ elements was defined by MacPherson in \cite{Mac} and is often
called the {\em MacPhersonian}.  The limit of the finite MacPhersonians gives an infinite
poset $\MacP(k,\infty)$, and its geometric realization
$\|\MacP(k,\infty)\|$ is the classifying space for rank $k$ matroid
bundles.  

Our main results are the Combinatorialization Theorem, which associates a matroid bundle to
a vector bundle, the Spherical
Quasifibration Theorem, which associates a spherical quasifibration to a matroid bundle, 
and the Comparison Theorem, which shows that the composition of these two
associations is essentially the forgetful functor.  More precisely, if $B$ is a regular cell
complex, and if $V_k(B)$, $M_k(B)$, and $Q_k(B)$ denote the isomorphism classes of rank $k$
vector bundles, matroid bundles, and spherical quasifibrations respectively, we construct maps
$$V_k(B) \to M_k(B) \to Q_k(B) 
$$
whose composite is the forgetful map given by deletion of the zero section.  The first map is
intimately related to the construction of a continuous map
$$\tilde{\mu} : G(k,\R^{\infty}) \to \|\MacP(k,\infty)\|
$$
and the Comparison Theorem leads to the following theorem.

\begin{A} 
The map 
$$\tilde{\mu}^* : H^*(\|\MacP(k,\infty)\|;\Z2) \to H^*(G(k,\R^{\infty});\Z2) 
$$
is a split surjection.
\end{A}

The mod 2 cohomology of $G(k,\R^\infty)$ is well-known: it is a polynomial algebra on the
Stiefel-Whitney classes $w_1,w_2, \dots , w_k$.  The above theorem gives Stiefel-Whitney
characteristic classes for matroid bundles.  

More generally, one can consider any rank $n$ oriented matroid $M^n$ as a combinatorial
analog to $\R^n$.  There is an associated {\em combinatorial Grassmannian}
$\Gamma(k,M^n)$ which is a partially ordered set of rank $k$ ``subspaces'' of
 $M^n$. The MacPhersonian arises in this way as $\MacP(k,n)=\Gamma(k, M_c)$, 
where $M_c$ is the unique rank $n$ oriented matroid  with elements
 $\{1,2,\ldots,n\}$.  
If $M^n$ is {\em realizable}, any realization induces a simplicial map
$$\tilde{\mu} : G(k,\Rn) \to \|\Gamma(k,M^n)\|
$$
from a triangulation of the real Grassmannian of $k$-planes in $\Rn$ to the 
geometric realization of the combinatorial Grassmannian, constructed by the 
same method as for the special case of the MacPhersonian.

\begin{B} 
The map 
$$\tilde{\mu}^* : H^*(\|\Gamma(k,M^n)\|;\Z2)\to H^*(G(k,\R^n);\Z2) 
$$
is a split surjection.
\end{B}

There is a natural combinatorial analog to an orientation of a real vector space leading to
the definition of an {\em oriented combinatorial Grassmannian} $\tilde{\Gamma}(k,M^n)$
analogous to the  Grassmannian $\tilde{G}(k,\Rn)$ of oriented $k$-planes in $\Rn$.  For any realizable oriented matroid 
$M^n$, there
is a combinatorialization map $\tilde{\mu}$ from $\tilde{G}(k,\Rn)$ to
$\|\tilde{\Gamma}(k,M^n)\|$. 

\begin{C} \label{C}
The map 
$$\tilde{\mu}^* : H^*(\|\tilde{\Gamma}(k,M^n)\|;\Zz)\to H^*(\tilde G(k,\R^n);\Zz) 
$$
has the Euler class in its image.
\end{C}

These results suggest substantial power for a combinatorial approach to
characteristic  classes via matroid bundles. 

The Comparison Theorem also leads to results on homotopy groups of the
combinatorial Grassmannian. We show that the second
homotopy group of the MacPhersonian is the same as that of
the corresponding Grassmannian. (Similar results for the zero and first
homotopy groups of $\MacP(k,n)$ were previously known.). In addition, we
get results on the homotopy groups of general combinatorial Grassmannians (which need not even be connected (\cite{MR93})).

\begin{D} \label{D}
Let $M^n$ be 
a realized rank $k$ oriented 
matroid on $n$ elements. 
Let $p$ be a point in the image of $\tilde\mu: G(k,\Rn)\to\|\Gamma(k,M^n)\|$.
\begin{enumerate} 

\item  $\pi_k (\|\Gamma(k,M^n)\|,p)$ has $\Zz$ as a subgroup when $k$ is even and $n \geq 2k$.
\item  $\pi_i (\|\Gamma(k,M^n)\|,p)$ has $\Z2$ as a subquotient when $i \equiv 1,2$  (mod 8), $n-k \geq i$, and $k \geq i$. 
\item  $\pi_{4m} (\|\Gamma(k,M^n)\|,p)$ has $\Zz_{a_m}$ as a subquotient when $m > 0$, $n-k \geq 4m$, and $k \geq 4m$.  Here $a_m$
is the denominator of $B_m/4m$ expressed as a fraction in lowest terms, and $B_m$ is the $m$-th Bernoulli number.
\end{enumerate}
\end{D}
Together with a previously known stability result (\cite{homotopy}), this
implies a similar statement about the infinite MacPhersonian. These are
the first results giving nontrivial $\pi_i(\|\Gamma(k,M)\|,p)$ for general
realizable $M$ for any
$i$.

Our results have potential interest in combinatorics and topology.  A major focus of work on
oriented matroids has been construction of oriented matroids behaving very differently from those in the image of $\tilde\mu$, 
but our partial computations of cohomology and homotopy groups should be viewed
as an attempt to tame the beast.  As far as geometric topology is concerned, a matroid bundle over a
finite cell complex is a purely finite gadget, and we have shown that matroid bundles have
characteristic class information.  MacPherson has conjectured that in fact all characteristic classes
of a vector bundle should be carried by the associated matroid bundle.  Also, the study of matroid
bundles is a necessary first step in the study of CD manifolds.

The paper is organized as follows.  We review the theory of oriented matroids and
develop the foundations of matroid bundles.  We show
that any vector bundle whose base space is a regular cell 
complex defines an isomorphism class of matroid bundles, thus
passing from topology to combinatorics.  We then construct a  combinatorial ``sphere
bundle'' (actually, a spherical quasifibration) associated to any matroid bundle, and thus
pass from combinatorics back to topology.  This construction is described in 
Section~\ref{constr}. This Combinatorialization Theorem depends on a deep result from
combinatorics (the Topological Realization Theorem) and a deep result from topology (Quillen's
Theorem B).  
From this construction we derive Stiefel-Whitney classes and Euler classes for matroid
bundles.  The Comparison Theorem (see Section \ref{vb and cvb}) is proven by constructing a map of
spherical quasifibrations from the canonical sphere bundle over the real Grassmannian to the
``sphere bundle'' over the combinatorial Grassmannian, which implies that the Stiefel-Whitney and
Euler classes  we defined for matroid bundles map under
$\tilde{\mu}^*$ to the analogous classes for real vector bundles.  This completes the
proof of the Theorems A, B, and C. Theorem D (Section \ref{sec:homotopy}) is a consequence of the Combinatorialization Theorem,
the Spherical Quasifibration Theorem, and the Comparison Theorem and homotopy theoretic results on the image of the
$J$-homomorphism.

We also show that these  combinatorial characteristic classes have
interpretations analogous to those of their classical counterparts, as obstructions to the
existence of  combinatorial ``orientations'' (Theorem~\ref{orient}) and
combinatorial ``independent sets of  vector fields'' (Section~\ref{obstr}).

We outline below the organization of the paper.
\begin{tabbing}
{1} Preliminaries\\
\qquad \={1.1} Oriented matroids\\
\>{1.2} Partially ordered sets\\
\>{1.3} Combinatorial Grassmannians\\
{2} Matroid bundles\\
\>{2.1} Matroid bundles and their morphisms\\
\>{2.2} Combinatorializing vector bundles: the Combinatorialization
Theorem\\
\>{2.3} Combinatorial sphere and disk bundles\\
{3} Combinatorial bundles are quasifibrations\\
\>{3.1} Quasifibrations\\
\>{3.2} The Spherical Quasifibration Theorem\\
{4} Stiefel-Whitney classes and Euler classes of matroid bundles\\
{5} Vector bundles vs. matroid bundles: the Comparison Theorem\\
{6} Homotopy groups of the combinatorial Grassmannian\\ 
{7} Vector fields and characteristic classes\\
{8} Some open questions\\
Appendix A Topological maps from combinatorial ones\\
Appendix B Babson's criterion
\end{tabbing}

\section{Preliminaries}

\subsection{Oriented Matroids}

For introductions to oriented matroids, see~\cite{BLSWZ} and \cite{Mac}. There 
are several equivalent axiomatizations of oriented matroids; throughout
this paper we will use the \emph{covector axioms} \cite[p.159]{BLSWZ}.

\begin{defn} An \textbf{oriented matroid} $M$ is a finite set $E(M)$ and a
subset $$\Vv^*(M) \subseteq \{-,0,+\}^{E(M)}$$  satisfying the following axioms:
\begin{enumerate}
\item $0\in\Vv^*(M)$.
\item If $X \in \Vv^*(M)$, then $-X \in \Vv^*(M)$.
\item {(\em Composition)} If $X, Y\in \Vv^*(M),$ then the function 
$$\begin{array}{lccl}
X\circ Y:& E(M)&\rightarrow&\{-,0,+\}\\
\\
         & e&\mapsto&\left\{\begin{array}{ll}
                                X(e)&\mbox{if }X(e)\neq 0\\
                                Y(e)&\mbox{otherwise}
                                \end{array}\right.
\end{array}$$
is in $\Vv^*(M)$.
\item {(\em Elimination)} If $X(e)=+$ and $Y(e)=-$, then there is a $Z\in \Vv^*(M)$ such that $Z(e) = 0$,
 and for all 
$f \in E(M)$ for which $X(f)$ and $Y(f)$ are not of opposite signs, $Z(f) = X \circ Y(f)$.

\end{enumerate}
$E(M)$ is called the set of \textbf{elements} of $M$.
$\Vv^*(M)$ is the set of \textbf{covectors} of $M$. 
\end{defn}

The motivating example: consider $n$ linear forms $\{\phi_1,\ldots,\phi_n\}$
on a finite dimensional real vector space $V$.  To any $p\in V$ we
associate a sign vector
$X(p)=(\sign\phi_1(p),\ldots,\sign\phi_n(p))\in\{-,0,+\}^n$. The set
$\{X(p):p\in V\}$ is the set of covectors of an oriented matroid $M$ with elements
$\{1,2,\dots,n \}$.  The set $\{\phi_1,\ldots,\phi_n\}$ is called a \textbf{realization} of
$M$. Any $M$ arising in this way is called \textbf{realizable}. Note that any of 
the $\phi_i$ can be multiplied by a positive scalar without changing $M$. Thus
we can represent a realizable oriented matroid by an arrangement of hyperplanes
$\{ \phi_i^{-1}(0) \}$ and a distinguished side $\phi^{-1}_i(\R^+)$
for each hyperplane. (If $\phi_i=0$ then the corresponding ``hyperplane'' is 
the \textbf{degenerate hyperplane} $V$ and the ``distinguished side'' is 
$\emptyset$.)  By considering a form as the inner product with a vector, we
arrive at yet another way of viewing an realizable oriented matroid: Take a finite
collection
$E =\{v_1,
\dots, v_n\}$ of vectors in a finite dimensional real inner product space $V$, then the
functions given by
$\{i
\mapsto
\text{sign}(v
\cdot v_i) : v
\in V\}$ are the covectors of an oriented matroid. 

\begin{defn} Let $M$ be an oriented matroid with elements $E$.  A subset $I$ of $E$ is
\textbf{independent} in $M$ if for every $e \in I$, there is a covector $X$ so that $X(e) \neq 0$,
but $X(I \backslash \{e\}) = 0$.  The \textbf{rank} of $M$ is the maximal order of a set of
independent elements of $M$. 
\end{defn}

If the oriented matroid arises from a set of vectors in a real inner product space, then
the rank of the oriented matroid equals the dimension of the span of $E$.

\begin{defn} \cite[pp. 133-134]{BLSWZ}
Let $A\subseteq E(M)$ where $M$ is an oriented matroid. Define two oriented
matroids whose elements are $E(M) \backslash A$.
\begin{enumerate}
\item  The covectors of the {\bf deletion $M \backslash A$} are 
$$
\{X|_{E(M)\backslash A}:
X\in\Vv^*(M) \}.
$$
\item  The covectors of the {\bf contraction $M /  A$} are 
$$
\{X|_{E(M)\backslash A }:
X\in\Vv^*(M) \text{ so that } X(a) = 0 \text{ for all } a \in A\}.
$$\end{enumerate}

\end{defn}

For a realizable oriented matroid, the deletion is realized by forgetting the linear
forms in $A$, while the contraction is realized by restricting the forms in $E(M) \backslash A$ to
the intersection of the zero sets of the forms in $A$.

Oriented matroids are connected to topology and geometry by the \emph{Topological
Representation Theorem} of Folkman and Lawrence.  This gives a way to associate a PL sphere to an oriented
matroid. We will state a weak version of the theorem here. First we take a definition from
\cite[\S5.1]{BLSWZ}.

\begin{defn} A \textbf{pseudosphere} $S$ in $S^{k-1}$ is the image of the equator $S^{k-2}$
under a homeomorphism $h : S^{k-1} \to S^{k-1}$. An \textbf{oriented pseudosphere} $S$ is a
pseudosphere together with a labeling $S^+$ and $S^-$ of the connected components of $S^{k-1}
\backslash S$.  Here $S^+$ and $S^-$ are called the \textbf{(open) sides of $S$}.  A
\textbf{signed arrangement of  pseudospheres} is a  finite multiset $\Aa=(S_e)_{e\in E}$, where for
each $e$, $S_e$ is  an oriented  pseudosphere of $S^{k-1}$, provided the following three
conditions hold:
\begin{enumerate}
\item  $S_A = \cap_{e \in A} S_e$ is homeomorphic to a sphere, for all subsets $A$ of $E$.
\item If $S_A \not\subseteq S_e$, for $A \subseteq E$ and $e \in E$, then $S_A \cap S_e$ is an
oriented pseudosphere in $S_A$ with sides $S_A \cap S_e^+$ and $S_A \cap S_e^-$.
\item The intersection of a collection of closed sides is either a sphere or a ball.
\end{enumerate}  $\Aa$ is \textbf{essential} if
$\cap_{e\in E} S_e=\emptyset$.

 \end{defn}

For instance, if $\{\phi_1,\ldots,\phi_n\}$ is a realization of $M$
in $\Rk$, then the arrangement of equators $\{\phi_i^{-1}(0)\cap S^{k-1}\}_{i \in \{1,2, \dots,
n\} }$ is a signed arrangement of pseudospheres. This arrangement is essential
if and only if the rank of $M$ is $k$.

Any essential signed arrangement of pseudospheres chops
$S^{k-1}$ into a regular cell complex.  The Topological Realization Theorem shows that the
cells of this complex give the nonzero covectors of an oriented matroid and that essentially all
oriented matroids arise in this way.

Let $\Aa = (S_e)_{e \in E}$ be a essential signed arrangement of psuedospheres in $S^{k-1}$. 
Define a function
$$
\sigma : S^{k-1} \to \{+,0,-\}^E
$$
by $\sigma(x)(e) = +, 0,$ or $-$
depending on whether $x \in S^+_e,S_e, \text{ or }
S^-_e$.  Then 
$$\{\sigma^{-1}(X) : X \in \sigma(S^{k-1})\}
$$ 
gives a regular cell decomposition of $S^{k-1}$ and  $\Ll(\Aa) =
\sigma(S^{k-1}) \cup \{0\}$ is the set of covectors of an oriented
matroid on $E$ (see \cite[\S 5.1]{BLSWZ}).

A \textbf{loop} in an oriented matroid is an element $e$ so that  $X(e) =
0$ for all covectors $X$. An oriented matroid is \textbf{loopfree} if it
has no loops.

\begin{TRP} (cf.~\cite{FL}, \cite{BLSWZ})
\begin{enumerate}
\item  If $M$ is a loopfree, rank $k$ oriented matroid on $E$, there is 
an essential signed arrangement
$\Aa=(S_e)_{e\in E}$ of pseudospheres in $S^{k-1}$ with  $\Ll(\Aa) =
\Vv^*(M)$.
\item If $\Aa=(S_e)_{e\in E}$ is an essential signed
arrangement of pseudospheres in $S^{k-1}$, then $\Ll(\Aa)$ is the  set of
covectors of a loopfree, rank $k$ oriented matroid on $E$.
\end{enumerate}
\end{TRP}

\subsection{Partially ordered sets}
Oriented matroids and geometric topology are related via partially
ordered sets, or \emph{posets}.  There are several partial  orders
associated to oriented matroids, where moving up in the partial order
corresponds to some  notion of moving into more general 
position.  The transition from posets to geometric topology is  through a
functor called \emph{geometric realization}.

\begin{defn} We define three partial orders:
\begin{enumerate}
\item On $\{-,0,+\}$:  The only strict inequalities are $+>0$ and $->0$.
\item On $\{-,0,+\}^E$:  Define $X \geq Y$ if $X(e) \geq Y(e)$ for all $e \in E$.
\item On the set of oriented matroids on a set $E$:  Define $M_1 \geq M_2$ if 
for every
$X\in\Vv^*(M_2)$ there is some $Y\in\Vv^*(M_1)$ such that $Y\geq X$.
\end{enumerate} 

\end{defn}

This last relation is sometimes described by saying 
$M_2$ is a \textbf{weak map image} of $M_1$, or that 
$M_2$ is a \textbf{specialization} of $M_1$. If
$\{\phi_1,\ldots,\phi_n\}$ is a  realization of a rank $k$ oriented matroid $M_1$ and
$\{\xi_1,\ldots,\xi_n\}$ is a  realization of $M_2$, then there is a weak map
$M_1\rightsquigarrow M_2$ if and only if $\sign(\phi_{i_1}\wedge\cdots\wedge\phi_{i_k}) \geq
\sign(\xi_{i_1}\wedge\cdots\wedge\xi_{i_k})$ for every $\{i_i,\ldots, i_k\}\subseteq
\{1,\dots ,n \}$.

An (abstract) \textbf{simplicial complex} $K$ is a collection of non-empty finite sets,
closed under proper inclusion.  The elements of $K$ are called \emph{simplices}.  An $i$-simplex is
an simplex with
$i+1$ elements, and
$K^i
\subset K$ is the set of $i$-simplices.  Note that a simplicial complex is a poset, with
partial order given by inclusion.
  A \textbf{regular cell
complex}
$B$ is a CW complex so that every cell
$e$ has a characteristic map $D^n \to \overline{e}$ which is a homeomorphism.  The
face lattice $\Ff(B)$ is the poset of  closed cells of $B$, ordered by inclusion. A simplicial
complex
$K$ determines a regular cell complex, which we  denote by $\|K\|$.  A \textbf{chain} in a poset $P$
is a non-empty, finite, totally ordered subset of $P$.   The
\textbf{order complex}
$\Delta P$ of a poset
$P$ is the simplicial complex  whose simplices are the chains in $P$.

Thus there are functors
$$\text{Posets} \xrightarrow{\Delta} \text{Simplicial Complexes} \xrightarrow{\|~\|}
\text{Regular Cell Complexes} \xrightarrow{\Ff} \text{Posets}
$$

For a poset $P$, we write $\|P\|$ for $\|\Delta P\|$, and call it
the \textbf{geometric realization} of $P$.  If $P^{\text{op}}$ is the poset given by reversing the inequalities,
then $\| P\| \cong \|P^{\text{op}}\|$.  For a regular cell complex $B$,  the simplicial complex given by
$\Delta
\Ff (B)$ is called the \textbf{barycentric subdivision} of $B$ and  there is a homeomorphism
$ B \cong  \| \Ff (B)\| $ under which every closed cell $\sigma$ of $B$ maps 
homeomorphically to the subcomplex  
$\|  \Ff
(B)_{\leq \sigma}\|$.  Here if $p$ is an element of a poset $P$, then $P_{\leq p} = \{q \in P : q \leq p \}$.
For a simplicial complex $K$, there is a {\em canonical} homeomorphism $ \|K\| \cong  \| \Ff (K)\|
$.

\begin{OHL}  If $f,g : P \to Q$ are poset maps so that for all $p \in P$, $f(p) \geq g(p)$ then
$\|f\|$ is homotopic to  $\| g\|$.
\end{OHL}

The proof is easy.  Simply let $(1 > 0)$ be the poset with elements $1$ and $0$, and the only
strict inequality being $1 > 0$. Then apply $\|~\|$ to the poset map $h : P \times (1 > 0)
\to Q$ defined by
$h(p,1) = f(p)$ and $h(p,0)= g(p)$.  A consequence of the lemma is that
if $P$ has a maximal or minimal element, then $\|P\|$ is contractible.

\subsection{Combinatorial Grassmannians}\label{defn:combgrass}

The \textbf{MacPhersonian} $\MacP(k,n)$ is the poset of rank $k$ oriented matroids on the set
$\{1,2,\ldots,n\}$, with $M_1 \geq M_2$ if there is a weak map $M_1 \rightsquigarrow    M_2$.    
There is an obvious embedding 
$\MacP(k,n)\hookrightarrow\MacP(k,n+1)$, by adding $n+1$ as a loop to each oriented matroid.  We
identify $\MacP(k,n)$ with its image under this embedding and define
$\MacP(k,\infty)$  to be the direct limit over $n$  of the $\MacP(k,n)$.
For any rank $n$ real vector space
$W$ with a fixed basis $\{w_1,\ldots,w_n\}$ (and therefore a
fixed inner product) there 
is a canonical function 
$$\mu : G(k,W) \to \MacP(k,n)
$$
given by intersecting each $V\in G(k,W)$ with the hyperplanes 
$\{w_1^\perp,\ldots, w_n^\perp\}$
 and considering the corresponding oriented matroid $\mu(V)$.  Equivalently one projects the basis
of $W$ onto $V$ and thus obtains an oriented matroid.
  This function $\mu$ is definitely not
continuous, but the point inverses give an interesting decomposition  of the
Grassmannian. (See~\cite[Section 2.4]{BLSWZ} for more about this decomposition.)

For future reference, we record a generalization of the MacPhersonian.

\begin{defn} There is a \textbf{strong map} $M_1\rightarrow M_2$ if 
$\Vv^*(M_2)\subseteq \Vv^*(M_1)$.
\end{defn}

For instance, if $\{\phi_1,\ldots,\phi_n\}$ is a 
realization of $M$ in $V$ and $W$ is a subspace of $V$, then 
$\{\phi_1|_W,\ldots,\phi_n|_W\}$ is the realization of a strong map 
image of $M$.

If
$M$ is an  oriented matroid, the
\textbf{combinatorial Grassmannian} $\Gg(k,M)$ is the poset of all rank $k$ strong map
images  of $M$, 
 with partial order given 
by weak maps.  If $M$ is the \emph{coordinate} rank $n$
oriented matroid, i.e., the unique rank $n$ oriented matroid with elements
$\{1,2,\ldots,n\}$, then $\Gg(k, M)$ is the MacPhersonian $\MacP(k,n)$. 

Let $M$ be a realizable rank $n$ oriented matroid, and let
\mbox{$\{v_1,\dots,v_m\} \subset \Rn$} realize  $M$.  Then, as above, there is a function  $\mu :
G(k,\Rn) \to \Gg(k,M)$ given by intersecting each vector space 
$V \in G(k,\Rn)$ with the oriented hyperplanes $\{v_1^\perp,\dots,v_m^\perp\}$ and
taking the corresponding oriented matroid $\mu(V)$.

\section{Matroid bundles}

We will define a combinatorial vector bundle over a regular cell complex $B$, to be an assignment of
an oriented matroid $M(e)$ to every cell $e$ of $B$ so that if $f$ is a face of $e$ then $M(e)$
weak maps to $M(f)$.   To see intuitively how a combinatorial vector bundle is derived from a vector
bundle over a finite regular cell complex $B$, note that if we fix a
metric for the bundle and a finite set $S$ of sections, then for every 
element $b$ of the base space, $\{s(b)\}_{s\in S}$ determines an oriented 
matroid $M_b$ with elements $S$. If we choose such an $S$ so that for every $b$,
$\{s(b)\}_{s\in S}$ spans the fiber over $b$ and if the function $b\mapsto M_b$ is
constant on the interior of each cell of $B$, then these oriented matroids 
determine a combinatorial vector bundle structure on $B$. In this case we
say $S$ is \emph{tame} with respect to $B$.
We will show that tame
sections exist (perhaps after a subdivision of the base space)
for vector bundles over finite-dimensional regular cell complexes.

\subsection{Matroid bundles and their morphisms}

The following is a generalization of the definition in~\cite{Mac}.

\begin{defn} A rank $k$ \textbf{matroid bundle} $\xi =
(B,{\Mm})$ is a poset $B$ and
a poset map ${\Mm}: B \to\MacP(k,\infty)$.
\end{defn}

\begin{defn}  The \textbf{universal rank $k$ matroid bundle} is
$$\gamma_k = (\MacP(k,\infty),\text{Id})$$
\end{defn}

\begin{defn}  A rank $k$
\textbf{combinatorial vector bundle} $\xi =  (B,{\cal M})$ is a piecewise-linear (PL)
cell complex $B$ and a poset map
${\cal M}$ from the set of cells of a PL subdivision of $B$, ordered by
inclusion, to $\MacP(k,\infty)$.  In other
words a combinatorial vector bundle $(B,\Mm)$ is a matroid bundle
$(\Ff(B'),\Mm)$ where $\Ff(B')$ is the poset of
cells of a PL subdivision $B'$ of $B$.
\end{defn}

Note every regular cell complex can be given the structure of a PL space via a
barycentric subdivision.
  
A matroid bundle $(B,\Mm)$ gives a combinatorial vector bundle $(\|\Delta
B\|,\Mm')$,
where $\Mm'(b_0<b_1<\cdots<b_m)=\Mm(b_m)$.
A combinatorial vector bundle $\xi = (B,\Mm : \Ff(B') \to
\MacP(k,\infty))$ induces a combinatorial vector bundle
$\xi'=(B,\Mm' : \Ff(B'') \to \MacP(k,\infty))$ for
 any PL subdivision
$B''$ of $B'$, by sending any cell $\sigma$ of $B''$ to
$\Mm(\delta(\sigma))$,
where $\delta(\sigma)$ is the smallest cell   of $B'$ containing
$\sigma$.
Two   combinatorial
vector bundles over a PL cell complex
$B$ are \textbf{equivalent} if they are equivalent
under the equivalence relation generated by PL
subdivision.  Two   matroid bundles over a poset $B$ are \textbf{equivalent} if the
associated combinatorial vector bundles are equivalent.  Clearly there is a bijective
correspondence between equivalence classes of matroid bundles over a poset and combinatorial
vector bundles over its geometric realization, and henceforth we blur the distinction.
                                 
\begin{defn} If $(B_1,\Mm_1)$ and $(B_2,\Mm_2)$ are two matroid bundles,
a \textbf{morphism} from $(B_1,\Mm_1)$ to  $(B_2,\Mm_2)$ is a triple
$(f,[C_f,\Mm_f])$, where
$f$ is a 
PL map  from $\| B_1\|$ to $\|B_2\|$, and
$[C_f,M_f]$ is
an equivalence class of combinatorial vector bundles over the  mapping
cylinder of $f$,
where $\Mm_f$  restricts to structures  equivalent to $(B_i,\Mm_i)$ at
either end.  Such a
morphism is a
\textbf{morphism covering $f$}.   A morphism covering the identity on
$\|B\|$ is called a
\textbf{$B$-isomorphism}. Clearly equivalent bundles over $B$ are $B$-isomorphic.
\end{defn}

\begin{defn} Let  $B_1$ be a poset, $\xi=(B_2,\Mm)$  a matroid bundle,
and $f:\|B_1\|\to \|B_2\|$  a PL map. Then $f$ is simplicial with respect to
some PL subdivisions $B_1'$ and $B_2'$ of $\|B_1\|$ and $\|B_2\|$. The composite
map $\Ff(B_1') \to \Ff(B_2') \to \MacP(k,\infty)$ defines a combinatorial vector bundle
over $\|B_1'\|$, and any subdivision of $B_1'$ gives an equivalent bundle. Thus $f$
defines an equivalence class of combinatorial vector bundles over $\|B_1\|$. We call any of the
corresponding matroid bundles a 
\textbf{pullback} of $\xi$ by $f$ and denote it $f^*(\xi)$.

\end{defn}

\begin{defn} If $(B_i,\Mm_i)$, $i\in\{1,2,3\}$ are rank $k$ matroid bundles and
$(f:B_1\to B_2,[C_f,\Mm_f])$ and $(g:B_2\to B_3,[C_g,\Mm_g])$ are
morphisms between them, then
consider the space $C_f\cup_{B_2}C_g$ obtained from the disjoint union of
$C_f$ and $C_g$ by identifying each $b\in B_2\subset C_f$ with $b\times
0\in C_g$. The matroid bundle structures on $C_f$ and $C_g$ define a matroid
bundle structure on this space. The pullback of the PL map
$c: C_{g\circ f}\to C_f\cup_{B_2}C_g$ defined by
$$c[b,i]=\left\{ \begin{array}{ll}
[b,2i]&\mbox{if $ i\leq 1/2$}\\
\mbox{$[f(b), 2(i-\frac{1}{2})]$}&\mbox{if $ i\geq 1/2$}
\end{array}\right.$$
and $c[b]=[b]$ for all $b\in B_3$, 
defines an equivalence class $[C_{g\circ f}, \Mm_{g\circ f}]$ of matroid
bundles which we call the  \textbf{composition} of the
two original morphisms.

\end{defn}
Note that two bundles are $B$-isomorphic if and only if there is a
combinatorial vector bundle over
$\|B\| \times I$ restricting on the ends to bundles equivalent to the
original ones.

The following familiar properties of bundles are easily verified.

\begin{prop} 
\begin{enumerate}
\item  Let $B_1$ and $B_2$ be posets.
There is a morphism $\xi_1  \to \xi_2$
covering a PL-map
$f :
\|B_1\| \to
\|B_2\|$ if and only if $f^*\xi_2$ is $B_1$-isomorphic to $\xi_1$.
\item  If $\phi,\psi : A \to B$ are poset maps so that $\|\phi\|$ and
$\|\psi\|$ are homotopic, and if $\xi$ is
a matroid bundle over $B$, then $\phi^*\xi$ is $A$-isomorphic to
$\psi^*\xi$.
\item  Every rank $k$ matroid bundle is a pullback of the universal rank
$k$ bundle.   
\item  If $\xi$ is a rank $k$ matroid bundle and there are two morphisms
$\xi \to
\gamma_k$ covering $f$ and $g$ respectively, then $f$ and $g$ are
homotopic.
\end{enumerate}
\end{prop}

 Recall that $[X,Y]$ is the set of homotopy classes of maps from $X$ to
$Y$.

\begin{corollary}  For a regular cell complex $B$, let $M_k(B)$ be the
set
of $B$-isomorphism classes
of rank $k$ combinatorial vector bundles over $B$.  Then
\begin{align*}
M_k(B) & \to  [B, \|\MacP(k,\infty)\|] \\
[\xi] & \mapsto [\|\Mm(\xi)\|]
\end{align*}
is a bijection, natural in $B$.  The inverse is defined by applying the
simplicial approximation theorem to [f] to obtain a subdivision $B'$ of $B$, a poset map $f' : B' \to
\Ff\|\MacP(k,\infty)\|$, and thus a matroid bundle $(B',\Mm)$ where $\Mm(\sigma)$ is the maximal vertex of
$f'(\sigma)$.
\end{corollary}
    
Thus the classifying space for rank $k$  combinatorial vector bundles is
$\|\MacP(k,\infty)\|$
with universal element $\gamma_k$.

Matroid bundles arise in combinatorics in a variety of ways:
see~\cite{AB}
for examples. In addition, any real vector bundle yields a combinatorial
vector bundle,    
as described in the following section.

\subsection{Combinatorializing vector bundles:
the Combinatorialization Theorem}\label{combinatorialized}

For any real rank
$k$  vector bundle $\xi =\mbox{$(p:E\rightarrow B)$}$ over a paracompact base space there
is a bundle map 
$$
\begin{CD}
E @>\tilde c>> E(k, \Rinfty)\\
@VVV @VVV \\
B @>c>> G(k, \Rinfty)
\end{CD} 
$$
to the canonical bundle over the Grassmannian of $k$-planes in $\Rinfty$,
 and $c$ is determined up to homotopy (cf.~\cite{MS} Ch. 5).  If $B$ is the underlying space of a regular cell
complex, we 
 call the map $c$ \textbf{tame} if  $\mu\circ c$ is constant on the interior of each cell. 
Such a tame classifying map gives a  combinatorial vector bundle $c(\xi) = (\Ff(B),\Mm)$ by defining $\Mm(\sigma)
=
\mu(c(\mbox{int }\sigma))$. 
Here $\Mm$ is a poset map since  if $\sigma$ is a face of $\tau$ , then 
$$c(\mbox{int }\sigma) \cap
\overline{c(\mbox{int }\tau)}
\neq
\emptyset$$ so $\mu(c(\mbox{int }\sigma)) \leq \mu(c(\mbox{int }\tau)) $ by
\cite[2.4.6]{BLSWZ}. A subdivision of the cell complex leads to
an equivalent combinatorial vector bundle.

Generalizing the notation a bit, for any vector space $V$, let $G(k,V)$ be the Grassmannian of 
$k$-planes in $V$.  A \textbf{classifying map} for a vector bundle $\xi = (p:E \to B)$ is a map $c : B
\to G(k,V)$ covered by a bundle map from $\xi$ to the canonical bundle.  For any finite set
$F$, let
$\MacP(k,F)$ be the poset of rank
$k$ oriented matroids with elements $F$.  For any set $A$, let $\MacP(k,A)$ be the direct limit of
$\MacP(k,F)$, taken over all finite subsets $F$ of $A$.  If $V$ is a vector space with a finite
basis $A$, define $\mu: G(k,V) \to \MacP(k, A)$ as we did in
Section~\ref{defn:combgrass}, while if the basis $A$ is infinite, define $\mu$ so
that it restricts to $\mu: G(k,\text{Span } F) \to \MacP(k, F)$ for all finite
subsets $F$ of $A$. If $V$ is a vector space, a \textbf{tame classifying map} is a
classifying map
$c: B \to G(k,V)$  so that  $\mu\circ c$ is
constant on the interior of each cell.

\begin{theorem}[Combinatorialization Theorem]\label{combinatorialize}
 Let  $\xi=(p: E\rightarrow B)$ be a rank
$k$ real  vector bundle, where $B$ is the underlying space of a regular cell complex.
\begin{enumerate}
\item For $i = 0,1$, let $c_i : B \to G(k,V_i)$  be a tame classifying map for $\xi$.  Then there is
a tame classifying map $h : B \times I \to G(k,V_0 \oplus V_1)$ for $\xi\times
I$, restricting to $c_i$ on $B \times
\{i\}$. 
\item If $B$ is finite dimensional, any classifying map $c: B \to G(k,V)$ is homotopic to a classifying
map which is tame with respect to some simplicial subdivision of the barycentric subdivision of $B$. 
\end{enumerate}
\end{theorem}

\begin{proof}  Note that for a vector bundle $\xi = (p : E  \to B)$, specifying a classifying map 
$c: B \to G(k,V)$ together with a covering map $\tilde c : E \to E(k,V)$ is equivalent to
specifying a map for
$\hat{c} : E \to V
$
which is linear and injective on each fiber.   (Here 
$c(b)=\hat c(p^{-1}b)$.) If $V$ has a basis $A$, then $c$ is tame if and only if the function
\begin{align*} B & \to \text{subsets of $\{+,-,0\}^A$} \\
b & \mapsto \{a \mapsto \text{sign}(\hat c(e) \cdot a)\}_{e\in p^{-1}b}
\end{align*}
is constant on the interior of cells.

We prove something slightly more general than (1).  Suppose that $V$ has a basis $A$ and that $c_0,
c_1 : B
\to G(k,V)$ are two tame classifying maps for $\xi$ with covering maps $\widehat c_0, \widehat
c_1$.  Suppose also that for every
$e\in E$ and for every $a \in A$,   $\widehat c_0(e) \cdot a$ and $\widehat c_1(e)
 \cdot a$ do not have opposite signs.   
Then 
 $$\widehat{h_t}(e) = (1-t) \widehat{c_0}(e) + t\widehat{{c_1}}(e)  
\qquad 0 \leq t\leq 1
$$
defines a classifying map $h : B \times I \to G(k, V)$ for the bundle 
$E \times I \to
B\times I$ which is tame with respect to the product cell structure on
$B\times I$ and hence a tame homotopy between  $c_0$ and 
$c_1$.  Applying this to $V = V_0 \oplus V_1$ gives (1).

For part (2), it suffices to prove it for the universal case.  We will find a triangulation of
$G(k,\R^n)$ so that the identity map is tame, i.e., $\mu : G(k,\R^n) \to \MacP(k,n)$ is constant on the
interior of simplices.  Then given a classifying map $c : B \to G(k,V)$ where $B$ has dimension $r$,
by the cellular approximation theorem and the Schubert cell decomposition of the Grassmannian, there
is a homotopic map
$c' : B
\to G(k,V')$ where
$V'
\subseteq V$ is a vector space spanned by
$k+r$ elements of the basis. Finally, apply the Simplicial 
Approximation Theorem to map from the
barycentric subdivision of $B$ to the tame triangulation of $G(k,V)$.

Thus the following lemma applied to the coordinate oriented matroid (where $\Gamma(k,M) =
\MacP(k,n)$) completes the proof of the combinatorialization theorem.

\begin{lemma} \label{Mac} (cf. \cite{Mac}) Let $M$ be a realizable rank $n$ oriented
matroid with a fixed realization. Then there is a semi-algebraic triangulation $T$ of $G(k,\mathbb
R^n)$ and a simplicial map with respect to its barycentric subdivision 
$$\tilde\mu: G(k,\mathbb R^n)\rightarrow \|\Gg(k, M)\|$$
 such that for every vertex $v$ in the barycentric subdivision, one has $\tilde \mu(v) = \mu(v)$.
Furthermore, the homotopy class of 
$\tilde{\mu}$ is independent of the choice of semi-algebraic triangulation.
\end{lemma}

\begin{proof}  This is an application of Appendix \ref{PL}, theorems on  existence
and uniqueness of semi-algebraic triangulations, and the fact that (cf. 2.4.6
in~\cite{BLSWZ}) $\mu : G(k,\R^n) \to \|\Gamma(k,M)\|$ is upper semi-continuous.

A key tool is \cite[{\em Semi-algebraic triangulation theorem}]{Hiro} which proves 
that for any finite partition $\{U_i\}_{i\in I}$
 of a  bounded, semi-algebraic set $S$ into
semi-algebraic sets there exists a semi-algebraic triangulation of $S$
such that each $U_i$ is a union of the interiors of simplices.  Furthermore,  by
\cite[2.4]{Hiro}, for any two such semi-algebraic triangulations, there is a
semi-algebraic triangulation which is a common refinement.  Thus $S$ is a
$PL$ 
space.

In the case at hand   $\{\mu^{-1}(N)\}_{N \in  \Gg(k, M)}$ is a semi-algebraic
partition of
$G(k,\Rn)$.  In the language of Corollary \ref{lemma:triang}, the corresponding triangulation 
refines the stratification  given by the upper semi-continuous map
$\mu$, so the result follows.
\end{proof}

\end{proof}

\begin{corollary}  Let $B$ be a finite dimensional regular cell complex.  Let $V_k(B)$ be the set of
$B$-isomorphism classes of  rank $k$ vector bundles over
$B$.  There is a ``combinatorialization map''
$$C : V_k(B) \to M_k(B),
$$
natural in $B$, defined by sending a vector bundle to the combinatorial vector bundle 
given by a tame
classifying map.
\end{corollary}

\begin{proof}[Proof of Corollary.]  Let $\xi$ be a $k$-dimensional vector
bundle over $B$.
\begin{itemize}
\item {\em Existence}:  The combinatorialization theorem shows there is a tame 
classifying map
$$c :  B \to G(k,V).
$$
Define $C[\xi]$ to be the corresponding combinatorial vector bundle $[c(\xi)]$.
\item {\em Uniqueness}:  
Suppose
$$c_i :  B \to G(k,V_i)  \qquad i = 0,1
$$
are two tame classifying maps.  Applying the first part of the 
combinatorialization theorem gives a tame classifying map
$$ h :  B \times I \to G(k,V_0 \oplus V_1),
$$
for $\xi\times I$ restricting to $c_0$ and $c_1$ at either end.  The 
resulting combinatorial vector bundle over $B \times I$ gives an $B$-isomorphism 
between $c_0(\xi)$ and
$c_1(\xi)$.
\item {\em Naturality}:  Clear.
\end{itemize}
\end{proof}

We wish to extend the map $V_k(B) \to M_k(B)$ to infinite dimensional complexes. 
The problem, as shown in~\cite{notPL}, is that
$G(k,\Rinfty)$ has no tame triangulation, i.e. no
triangulation where $\mu$ is constant on simplices.  But we do have the following 
theorem.

\begin{theorem}  \label{mu} There is a map $\tilde \mu : G(k,\Rinfty) \to 
\|\MacP(k,\infty)\|$ 
which restricts to a map $G(k,\R^n) \to \|\MacP(k,n)\|$  given
by Lemma
\ref{Mac} for all $n$.  The homotopy class of $\tilde \mu$ is well-defined.

\begin{proof} This follows from Appendix \ref{PL} and results on existence of
semi-algebraic triangulations, by the same argument as the proof of
Lemma~\ref{Mac}.
\end{proof}

\end{theorem}

\begin{corollary}\label{cor:comb}
  Let $B$ be a regular cell complex.  There is a ``combinatorialization map"
$$C : V_k(B) \to M_k(B),$$
natural in $B$, which for finite-dimensional $B$ coincides with the map given by
sending a vector bundle to the combinatorial vector bundle given by a tame
classifying map.

\begin{proof} By replacing $B$ by $\|\Ff(B)\|$, we may assume that $B$ is the geometric realization
of a  simplicial complex.  Let  $c : B \to G(k, \Rinfty)$ be a classifying map for $\xi$. Apply the
simplicial approximation theorem to $\tilde \mu \circ c$  to find a subdivision $B'$ of $B$ and a
map $c' : \Ff (B') \to \MacP(k,\infty)$, where $\|c'\|$ is homotopic to $\tilde \mu \circ c$. 
Then set $C[\xi] = [c']$.
\end{proof}

\end{corollary}

\subsection{Combinatorial sphere and disk bundles}\label{constr}

Gelfand and MacPherson, in their combinatorial formula for the Pontrjagin
classes of a differentiable manifold, constructed a ``combinatorial
sphere bundle'' associated to a matroid bundle: 

\begin{defn}  For a matroid bundle $\xi=(B,\Mm)$, define posets
\begin{align*}
E(\xi) &=\{(\sigma, X): \sigma \in B, ~
        X\in \Vv^*(\Mm(\sigma))~\} \\
E_0(\xi) &=\{(\sigma, X): \sigma \in B, ~
        X\in \Vv^*(\Mm(\sigma))\backslash\{0\}~\}
\end{align*}
with $(\sigma, X) \geq (\sigma ', X')$ if $\sigma \geq \sigma'$ and $X \geq X'$.

The projection map $\pi_0:E_0(\xi)\rightarrow B$ and $\pi:E(\xi)\rightarrow B$
are the
\textbf{combinatorial sphere bundle} and \textbf{combinatorial disk bundle} 
associated to
$\xi$.
 \end{defn}

\begin{example}  {\em Warning: The geometric realization of the combinatorial 
sphere bundle
may not be a topological sphere bundle!}  The Topological Representation Theorem 
promises
that the realization of each fiber of $\pi_0$ over a vertex is a PL sphere. But in general, the
realization of $\pi_0$ is \emph{not} a topological sphere bundle, as we can see from the 
 example in Figure 1.

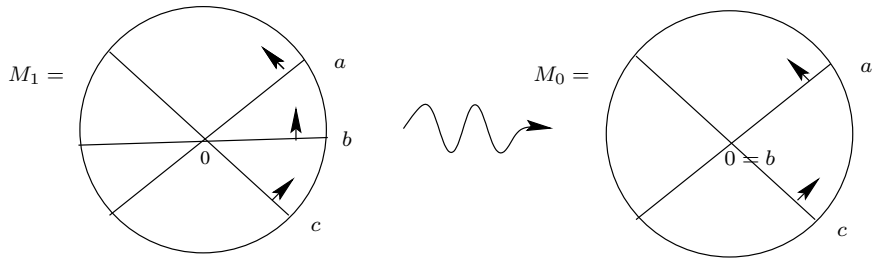
\begin{figure}[h]
\input{notbundle.pstex_t}
\caption{A sphere bundle need not be a sphere bundle.}
\end{figure}
This figure shows a weak map $M_1\rightsquigarrow M_0$
of rank 2 (realizable) oriented matroids with elements $\{a,b,c\}$; in
the second oriented matroid the element $b$ is the degenerate hyperplane
$0^\perp$. For each oriented matroid, the
nonzero covectors are given by the cell decomposition of the unit circle.
We can define a rank 2 matroid bundle over the poset $(1 > 0)$  by sending 1 to $M_1$ and 0
to $M_0$.  The total space of the associated sphere bundle will contain a 3-simplex
$\{(1, a^-b^+c^+),(1, a^-b^0c^+),(0, a^-b^0c^+),
(0, a^0b^0c^+)\}$. Hence the geometric realization of the combinatorial sphere bundle is not a 
topological circle bundle
over the 1-dimensional cell $\|1>0\|$.
\end{example}

It seems unlikely that the geometric realization of $\pi_0$ always gives a
fibration.  In  Section \ref{bundles_are_quasifibrations} we will show that it is
the next best thing, a spherical quasifibration.

\section{Combinatorial bundles are quasifibrations}\label{bundles_are_quasifibrations}

In this section we review the notion (due to Dold-Thom \cite{DT}) of a quasifibration and 
show that the combinatorial sphere bundle associated to a matroid bundle is a spherical
quasifibration.  A key tool is a criterion for the geometric realization of a poset map to be a
quasifibration.  This criterion was formulated in the Ph.D. thesis
\cite{Ba} of Eric Babson.  We give a proof of Babson's criterion in 
Appendix~\ref{eric}.
It is an application of Quillen's work on the foundations of algebraic
$K$-theory. 

\subsection{Quasifibrations}

\begin{defn}  A map $p : E \to B$ is a \textbf{fibration} if it has the homotopy lifting
property \cite{Whitehead}.  A map
$p : E
\to B$ is a
\textbf{quasifibration} if 
$$p_* : \pi_i(E,p^{-1}b,e) \to \pi_i(B,\{b\},b)$$
is an isomorphism for all $i \geq 0$, for all $b \in B$, and for all $e \in p^{-1}b$. 
 \end{defn}

\begin{defn}  A \textbf{spherical (quasi)-fibration of rank $k$} is a 
(quasi)-fibration $p_0 : E_0 \to B$ so that for all $b \in B$, $p_0^{-1}b$ has the weak homotopy 
type of $S^{k-1}$  (i.e., there is a map $S^{k-1} \to p_0^{-1}b$ inducing an isomorphism on
homotopy groups).
\end{defn}

A fibration
is a quasifibration.  An example of a quasifibration which is not a fibration is given by
collapsing a closed subinterval of an interval to a point. A (quasi)-fibration has a long exact
sequence in homotopy.  

A construction of Bourbaki \cite[\S I.7]{Whitehead} shows that every continuous map $f : E
\to B$ has the homotopy type of a fibration, i.e. there is a fibration $\pi_f : P_f \to B$ and
a homotopy equivalence $h: E \to P_f$ so that $f = \pi_f \circ h$.  Here 
$$P_f = \{(e,\alpha) \in E \times B^I : f(e) = \alpha(0)\}$$
$\pi_f(e,\alpha) = \alpha(1)$ and $h(e) = (e, \hbox{const}_e)$. For $b \in B$, $f^{-1}b$ is the
{\bf fiber above
$b$} and $\pi^{-1}_fb$ is the {\bf homotopy fiber above $p$}.  The homotopy long exact sequence of a
fibration and the five lemma give the following alternative (and perhaps better) definition of a
quasifibration.

\begin{prop}  A map $p: E \to B$ is a quasifibration if and only if  for all
$b \in B$
$$p^{-1}(b) \to {\pi_p}^{-1}(b)
$$
is a weak homotopy equivalence.
\end{prop}

\begin{defn}  A \textbf{morphism of
(quasi)-fibrations}
$p$ and $p'$ is a map $f: B \to B'$ and a (quasi)-fibration $E_f \to C_f$ over the
mapping cylinder of $f$ which restricts on the ends to $f$ and $f'$.
Such a morphism is called a \textbf{morphism covering $f$}.  A morphism covering the identity on $B$
is called a \textbf{$B$-isomorphism}.  A \textbf{CW-(quasi)-fibration}, respectively a
\textbf{$B$-CW-isomorphism}, is a (quasi)-fibration, respectively $B$-isomorphism, in
which the  domain of each (quasi)-fibration has the homotopy type of a
CW-complex.
\end{defn}

Let $p' : E' \to B$ and $p : E \to B$ be two maps.  A map $g: E' \to E$   is
\textbf{fiber-preserving} if $p \circ g = p'$.   A \textbf{fiber-preserving homotopy equivalence}
(f.p.h.e) is a homotopy equivalence $g: E' \to E$ which is fiber-preserving.  There
is also the notion of a \textbf{fiber-preserving weak homotopy equivalence} (f.p.w.h.e). 
Two maps
$g,g' : E' \to E$ are \textbf{fiberwise homotopic} if there is a homotopy
$G : E'
\times I
\to E$ between them so that for all $t$, $G(-,t)$ is fiber-preserving.  A  
fiber-preserving map $g : E'
\to E$ is a  \textbf{fiber homotopy equivalence} (f.h.e) if there is a  fiber-preserving 
map $h : E \to E'$  so that $g \circ h$ and $h \circ g$ are both fiberwise
homotopic to the identity.  One says that $p$ and $p'$ have the same  fiber
homotopy type.  Of course a fiber homotopy equivalence is a fiber-preserving
homotopy equivalence.

Two quasifibrations $p : E \to B$ and $p' : E' \to B$ are $B$-isomorphic  if and
only if they are equivalent under the equivalence relation generated by f.p.w.h.e.
Indeed if
$g : E'
\to E$ is a f.p.w.h.e., then the natural map $C_g \to B \times I$ shows
that
$p$ and
$p'$ are
$B$-isomorphic quasifibrations.  Conversely given a $B$-isomorphism $p'' : E'' \to B \times I$, the
inclusion map gives a f.p.w.h.e. from $p$ (or $p'$) to \mbox{$\text{pr}_B \circ p'' : E'' \to B$}.

Two fibrations $p: E \to B$ and $p' : E' \to B$ are $B$-isomorphic if and only if there is a fiber
homotopy equivalence between them, which occurs if and only if there is a fiber preserving homotopy
equivalence between them.  All of this is an elementary, if somewhat confusing, exercise in the
homotopy lifting property.  For a reference that $B$-isomorphism implies f.h.e., see
\cite[Theorem 7.25]{Whitehead} and for a reference that f.p.h.e. implies h.e., see \cite[Theorem 6.1]{D}.

The following theorem
shows that in terms of homotopy theory, 
there is really not much difference between fibrations and
quasifibrations.  This theorem, in slightly different language, is due to Stasheff \cite{S}.

\begin{theorem}\label{bijection}
  For a CW-complex $B$, let $Q(B)$ (respectively $F(B)$) be the set of $B$-CW-isomorphism
classes of CW-quasifibrations (respectively CW-fibrations) over $B$.  There is a bijection,
$$Q(B) \to F(B),
$$
given by converting a quasifibration $p$ into a fibration $\pi_p$.  The inverse is the
forgetful map, given by considering a fibration as a quasifibration.
\end{theorem}

\begin{proof}  This conversion process has two nice properties.  The first is that it
sends a fiber-preserving homotopy equivalence to a fiber-preserving homotopy equivalence.  The
second is that given a map
$p: E
\to B$ where
$E$ and $B$ have the homotopy type of a CW-complex, then $P_p$ also has the homotopy type of a
CW-complex \cite{M}.

We  need to see that the map $Q(B) \to F(B)$ is well-defined.  If $p'' : E'' \to B\times I$
is a $B$-isomorphism between quasifibrations $p: E \to B$ and $p': E' \to B$, then as above, there
is a a f.p.w.h.e. from $p$ (or $p'$) to 
$\text{pr}_B \circ p''$, which is a f.p.h.e by the CW-assumption.  Now this conversion process
takes a f.p.h.e to a f.p.h.e, and hence the corresponding fibrations are $B$-CW-isomorphic.

If one first converts and then forgets, one obtains a quasifibration which is f.p.h.e. to the
original one, and hence equivalent.  Conversely, if one has a fibration and converts it, the result
is a fibration equivalent to the original one.
\end{proof}

\begin{remark}  By the pullback, $F(-)$ is a contravariant functor from topological spaces
to sets.  However, since the pullback of a quasifibration need not be a quasifibration, it
is not clear that $Q(-)$ is a functor.  However, using the equivalence in the above
theorem, $Q(-)$ does give a functor from the category of spaces having the homotopy type of
CW-complexes to sets.  Furthermore, if the pullback of a  quasifibration $p : E  \to
B$ under a map $f : B' \to B$ happens to be a quasifibration, then the pullback
$f^*E \to B'$ represents the correct induced element of $Q(B')$. 
\end{remark}

\subsection{The spherical quasifibration theorem}

\begin{SQF}  \label{thm:quasi} For any matroid bundle $\xi = (B, \Mm)$,  the
geometric realizations of the combinatorial sphere and disk bundles
\begin{align*}  \|\pi_0\| & : \|E_0(\xi)\| \to \|B\| \\
\|\pi\| & : \|E(\xi)\| \to \|B\|
\end{align*}
are quasifibrations.
\end{SQF}
We use the following criterion for the geometric realization of a poset map  to be
a quasifibration.  

\begin{BC}  
If $f : P \to Q$ is a poset map satisfying both of the conditions below, then
$\|f\|$ is a quasifibration.
\begin{enumerate}
\item $\|f^{-1}q \cap P_{\leq
p}\|$ is contractible whenever $p\in P$, $q \in Q$, and $q \leq f(p)$.
\item $\|f^{-1}q \cap P_{\geq
p}\|$ is contractible whenever $p\in P$, $q \in Q$, and $q \geq f(p)$.
\end{enumerate}
\end{BC}

This 
criterion was formulated in the Ph.D. thesis \cite{Ba} of Eric Babson, and is an
application of Quillen's work on the foundations of algebraic $K$-theory.  We give a proof 
in Appendix \ref{eric}.  In this section we verify that the combinatorial bundles
satisfy Babson's criterion.

\begin{lemma}\label{contract}
 If $M\rightsquigarrow M'$ and $e$ is a nonzero element of $M'$
then $M/e\rightsquigarrow M'/e$.
\end{lemma}

\begin{proof}  Let $X' \in \Vv^*(M'/e) = \{Z' \in \Vv^*(M')~:~Z'(e) = 0\}$.  Since $e$ is nonzero in
$M'$, there are covectors $Z_1'$ and $Z_2'$ of $M'$ so that $Z_1'(e) = +$ and $Z_2'(e) = -$.  Since
$M\rightsquigarrow M'$, there are covectors $X_1$ and $X_2$ of $M$ so that $X_1 \geq X' \circ Z_1'
$ and $X_2
\geq X' \circ Z_2' $.  Since $X_1(e) = +$ and $X_2(e) = -$, we can apply the elimination
axiom in the definition of an oriented matroid.  What results is a covector $X$ of $M /e$
so that $X \geq X'$.
\end{proof}

\begin{lemma}  \label{lemma:babson2}
If $M\rightsquigarrow M'$,  $\text{\rm rank}~ M = \text{\rm rank} ~M'$, and $X$ 
is a nonzero covector of $M$,
then there is a nonzero covector $X'$ of $M'$ so that $X \geq X'$.
\end{lemma}

\begin{proof}
We induct on $\rank(M)$ and on the number of elements of $M$. When $\rank(M) = 1$, the existence of
$X'$ is easy.

If $\rank(M)>1$, it suffices to consider the case when $X$ is not maximal, since if $X$
is maximal, there is a nonzero covector $\overline{X}$ of $M$ so that $X > \overline{X}$ and
we replace $X$ by $\overline{X}$.  So assume that there is some nonzero
element $e$ of
$M$  such that $X(e)= 0$. We have two cases:
\begin{itemize}
\item If  $e$  is zero in $M'$, then $M\backslash e
\rightsquigarrow M'\backslash e$. Then by induction on the
number of elements we get $X'\in \Vv^*(M'\backslash e) = \Vv^*(M')$ such that $X\geq X'$.

\item If  $e$  is nonzero in $M'$, then by Lemma~\ref{contract}
we have $M/e\rightsquigarrow M'/e$. Since $X\in \Vv^*(M/e)$, induction on rank
gives a nonzero $X'\in \Vv^*(M'/e)\subset \Vv^*(M')$ such that $X\geq X'$.
\end{itemize}
\end{proof}

\begin{remark}  A consequence of this lemma is that if $M\rightsquigarrow M'$  and
$\text{rank}~M =
\text{rank}~M'$, there is a poset map $\Phi: \Vv^*(M) \to \Vv^*(M')$ which  maps
nonzero covectors to nonzero covectors and so that $X  \geq \Phi(X)$ for all $X$. 
Indeed, $\Phi(X)$ is defined to be the composition of all nonzero covectors of
$M'$ which are less than or equal to $X$.  This map $\Phi$ lends credence to the
intuition that a weak map corresponds to moving into special position. This map
and variations are explored further in~\cite{weakmaps}.
\end{remark}

\begin{lemma} Let $K$ be a simplicial decomposition of a compact PL  manifold
with boundary,  let $K^0 =\{\sigma \in K : \|\sigma\|\cap \|\partial K\|= 
\emptyset\}$, and  assume that $\partial K$ is \textbf{full}, i.e., any simplex
whose faces are all contained in
$\partial K$ is itself contained in $\partial K$.  Then
$\|K\| \simeq \|K^0\|$.
\end{lemma}

\begin{proof}  Enumerate the
simplices $\sigma_1, \sigma_2, \dots , \sigma_k$ of $\partial K$ so that the dimension is monotone
decreasing.  Let 
$$
K_i = K \backslash (K_{\geq\sigma_1} \cup \dots \cup K_{\geq\sigma_{i-1}}).
$$
Note $K_1 = K$ and $K_{k+1} = K^0$. 

 We show $\|K^0\|$ can be obtained from $\|K\|$ by a sequence of
elementary collapses and expansions. Define an elementary collapse  to be
{\bf inward} if it collapses out a pair of simplices $\omega$ and
$\omega\cup\{x\}$ with $\{x\}\not\in\partial K$. (We will use the term inward to apply to collapses of any complex, not just $K$.)  We show by induction on dimension and induction on $i$ that $\|K_{i+1}\|$ can be obtained from
$\|K_i\|$ by a sequence of elementary collapses and expansions. Both initial cases hold
because
$\partial K$ is full.

We will use $\sim$ to denote equivalent via a sequence of elementary collapses inwards and
elementary expansions.  
\begin{align*}
\|\link_{K_i} \sigma_i\| & \sim \|\link_{K} \sigma_i\| &\text{(induction on $i$ and $K \sim K_i$)}\\
& \sim \|\link_{K^0} \sigma_i\| &\text{(induction on dimension and $\|\link_{K} \sigma_i\|$ is a PL-ball)}\\
& \sim * &\text{(contractible by the last step and  $\sim *$ since $K^0 \cap \partial K = \emptyset$)}
\end{align*}

We leave for the reader to verify
that such a sequence of collapses inward and expansions from
$\|\link_{K_i}\sigma_i\|$ to a point gives a sequence of collapses inward and 
expansions from
$\|\star_{K_i}\sigma_i\|=\sigma_i\ast \|\link_{K_i}\sigma_i\|$ to 
$\|\link_{K_i}\sigma_i\|$, and hence from $\|K_i\|$ to $\|K_{i+1}\|$.

\end{proof}

\begin{proof}[Proof of Spherical Quasifibration Theorem]
 We apply Babson's Criterion, first with $P=E_0(\xi)$ and $Q=B$, then with
 $P=E(\xi)$ and $Q=B$. In the first case, 
let $(M,X)\in E_0(\xi)$ and $M'\in B$. Then $X \in
\VM
\backslash \{0\}$ and $\pi^{-1}(M')\cong\Vv^*(M')\backslash \{0\}$.

If $M\rightsquigarrow M'$, then $\pi^{-1}(M')\cap E_0(\xi)_{\leq (M,X)}$
is isomorphic to the poset of all covectors $X'$ of $M'$ such that $X\geq X'$.
This is the poset of all covectors in $M'$ corresponding to cells of
$$B_{M'}^X = (\bigcap_{e \in X^{-1}(+)
} \overline{S_e^+}) \cap (\bigcap_{e \in X^{-1}(-)
} \overline{S_e^-}) \cap
(\bigcap_{e\in X^{-1}(0)}  {S_e})
$$
given as a subcomplex of the pseudosphere picture of $M'$. By the Topological
Representation Theorem, this is either empty, a $PL$-sphere, or a
$PL$-ball.  It can't be a sphere since $X \neq 0$ and it is non-empty  by Lemma
\ref{lemma:babson2}.
 Thus the first condition of Babson's Criterion is fulfilled.

If $M'\rightsquigarrow M$, then $\pi^{-1}(M')\cap E_0(\xi)_{\geq (M,X)}$
is isomorphic to the poset of all covectors $X'$ of $M'$ such that $X'\geq X$.
This is the poset of all covectors in $M'$ corresponding to cells   in the
interior of the cell complex
$$B^{M'}_X = (\bigcap_{e \in X^{-1}(+)} \overline{S_e^+}) \cap (\bigcap_{e \in X^{-1}(-)}
\overline{S_e^-})
$$
given as a subcomplex of the pseudosphere picture of $M'$.  Now $B^{M'}_X$ must 
be empty or a $PL$-ball, but is in fact a
$PL$-ball of full rank since this is a non-empty intersection (containing  
cells corresponding to
covectors $X' \geq X > 0$).  
By the previous lemma, the realization of the poset of cells in the 
interior of this ball is contractible.
Thus the second condition of Babson's Criterion is fulfilled, and so the
realization of the combinatorial sphere bundle is a quasifibration.

In the case of the disk bundle, the first condition to check is trivial, since
the poset in question will have a unique minimal element. The second 
condition follows immediately from the proof of the second condition for sphere 
bundles.
\end{proof}

\begin{corollary}  Let $B$ be a regular cell complex and $Q_k(B)$ be the set of
$B$-isomorphism classes of rank $k$ spherical quasifibrations over $B$.  The geometric
realization of the combinatorial sphere bundle gives a well-defined map
$$\|E_0\|: M_k(B) \to Q_k(B)
$$
natural in $B$.
\end{corollary}

We now have two maps $V_k(B) \to Q_k(B)$, the map
above and the map given by deleting the zero section of a vector bundle. 
In Section \ref{vb and cvb} we will show they coincide.

\section{Stiefel-Whitney classes and Euler classes of matroid bundles}

Recall the axioms for Stiefel-Whitney classes \cite[\S 4]{MS}:
\begin{enumerate}
\item  For any vector bundle $\xi = ( p: E \to B)$ there are classes $w_i(\xi) \in
H^i(B;\Z2)$, with
$w_0(\xi) = 1$ and
$w_i(\xi) = 0$ when $i$ is larger than the fiber dimension.
\item  If $f: B' \to B$ is covered by a bundle map $\xi' \to \xi$, then 
$w_i(\xi') = f^*w_i(\xi)$.
\item  $w_n(\xi_0 \oplus \xi_1) = \sum_{i=0}^n w_i(\xi_0) \cup w_{n-i}(\xi_1).
$
\item  The first Stiefel-Whitney class of the canonical line bundle over $\R P^\infty$ is
non-trivial. 
\end{enumerate}

 The construction of the Stiefel-Whitney classes of a vector bundle
\mbox{$\xi =(p: E \to B)$} with fiber $\R^{k}$ 
given in \cite[\S 8]{MS} is 
$$w_i(\xi) = \phi^{-1}{\rm Sq}^i\phi(1) \in H^i(B;\Z2)
$$
where ${\rm Sq}^i$ is the $i$-th Steenrod square and 
$$\phi : H^*(B;\Z2) \to H^{*+k}(E,E_0;\Z2)
$$
is the Thom isomorphism.  We next review the construction of Stiefel-Whitney classes and Euler 
classes for spherical (quasi)-fibrations.

\begin{TI}  Let $p_0 :E_0 \to B$ be a rank $k$ spherical quasifibration.
Let $p : E \to B$ be a quasifibration with contractible fiber and a fiber-preserving
embedding $E_0 \to E$.  
\begin{enumerate}
\item  There is a class $U \in H^{k}(E,E_0; \Z2)$, so that for all $b\in B$,
$\rm{inc}^*U \in H^{k}(p^{-1}b,p_0^{-1}b; \Z2) \cong \Z2$ is non-zero.  Furthermore
\begin{align*}
\phi : H^i(B;\Z2) & \to H^{i+k}(E,E_0;\Z2) \\
\alpha & \mapsto p^*\alpha \cup U
\end{align*}
is an isomorphism for all $i$.
\item   If there is a class $U \in H^{k}(E,E_0)$, so that for all $b\in B$,
$\rm{inc}^*U \in H^{k}(p^{-1}b,p_0^{-1}b) \cong \Zz$ is a generator, then
\begin{align*}
\phi : H^i(B) & \to H^{i+k}(E,E_0) \\
\alpha & \mapsto p^*\alpha \cup U
\end{align*}
is an isomorphism for all $i$.
\end{enumerate}
\end{TI}

\begin{proof}  For any quasifibration $f : X \to Y$ and point $y \in Y$, there is a 
Serre spectral sequence 
$$ E_2^{i,j} = H^i(Y; H^j(f^{-1}y)) \Longrightarrow H^{i+j}X
$$
given by the Serre
spectral sequence of the associated fibration $\pi_f$.  If $f$ were  a fibration
to begin with, there are, a priori, two different Serre spectral sequences, since
$f$ can be considered as a quasifibration or as a fibration.  They coincide, since
if $f$ is a fibration then $f$ and $\pi_f$ have the same fiber homotopy type.

  The collapsing
of the relative Serre spectral sequence 
$$E_2^{i,j} = H^i(B;H^j(p^{-1}b,p_0^{-1}b;\Z2) \Longrightarrow H^{i+j}(E,E_0;\Z2)
$$
gives the Thom isomorphism, and the Thom class $U$ is the image of 1 under the Thom
isomorphism.

With integer coefficients, the same argument applies except that the 
$E^2$-term might have twisted coefficients.  However the existence of an integral Thom class 
in (2) guarantees that the coefficients are untwisted (look at $E^{0,k}_2$).
\end{proof}

\begin{defn}  $U$ is called the \textbf{Thom class} and $\phi$ is called the 
\textbf{Thom isomorphism}.  In case 2 above, the spherical (quasi)-fibration is called
\textbf{orientable} and a 
choice of Thom class
$U
\in H^{k}(E,E_0)$ is called an \textbf{orientation}.
\end{defn}

\begin{defn}  Suppose $\xi = (p_0 : E_0 \to B)$ is either a vector bundle  with
the $0$-section deleted, a combinatorial sphere bundle, or a spherical
(quasi)-fibration.  In the three cases respectively, let $p : E \to B$ be the
vector bundle, the combinatorial disk bundle, or the obvious map  $p:E
\to B$ from the mapping cylinder $E$ of $p_0$.  Then the \textbf{$i$-th
Stiefel-Whitney class of
$\xi$} is
$$w_i(\xi) = \phi^{-1}Sq^i\phi(1) \in H^i(B;\Z2).
$$
If $p_0$ is oriented with Thom class $U \in H^{k}(E,E_0)$, the \textbf{Euler class}
$$e(\xi) \in H^{k}(B)
$$
is the image of the Thom class under
$$H^{k}(E,E_0) \to H^{k}E \cong H^{k}B.
$$
\end{defn}

We next wish to show that Stiefel-Whitney classes and Euler classes  satisfy the
axioms and the usual properties, but first we had better make clear what is meant
by Whitney sum.  

\begin{defn}  If $\xi_1 = (B, \Mm_1 : B \to \MacP(i,E_1))$ and  $\xi_2 = (B, \Mm_2
: B \to
\MacP(j,E_2))$ are two matroid bundles, then the
\textbf{Whitney sum
$\xi_1
\oplus
\xi_2$} is the matroid bundle $(B, \Mm_1 \oplus \Mm_2 : B \to \MacP(i+j, E_1
\amalg E_2))$ sending each $b$ to the direct sum $\Mm_1(b)\oplus\Mm_2(b)$.  If
$\xi_1 = (p_1 :E_1 \to B)$ and $\xi_2 = (p_2 :E_2 \to B)$ are two spherical
(quasi)-fibrations then the \textbf{Whitney sum} is 
$$\xi_1 \oplus \xi_2 = (p_1 *_B p_2 : E_1 *_B E_2 \to B),
$$ 
where
$$E_1 *_B E_2 = \{[e_0,e_1,t] \in E_1 * E_2 : p_0(e_0) = p_1(e_1) \text{ whenever }
t \neq 0,1 \}
$$
is the 
fiberwise join.
\end{defn}

It is not difficult to show that the geometric realization of the combinatorial sphere
bundle of a Whitney sum of matroid bundles is the Whitney sum of the resulting spherical
quasifibrations.

\begin{prop}  The four axioms for Stiefel-Whitney classes are satisfied for vector
bundles, matroid bundles, and for spherical (quasi)-fibrations.
\end{prop}

\begin{proof} It suffices to prove the axioms hold for spherical fibrations. Axiom 1 holds
since  $Sq^0 = \text{Id}$ and  $Sq^i$ is zero on $H^{k}$ for $i > k$.  Axiom 2 is clear by
construction.

The Whitney sum formula (Axiom 3) holds since the Thom class $\phi(1)$ for the fiberwise join is the
external product of the Thom classes of the two summands and there is a sum formula for
the Steenrod squares.

As for Axiom 4, one may compute $w_1$ by restricting the canonical line bundle to
the circle.  Here the bundle is the M\"obius strip,
 which has non-trivial $w_1$ by
direct computation.
\end{proof}

\begin{remark}  While the axioms characterize the Stiefel-Whitney classes of vector bundles 
(due to the splitting principle), there is no assertion that the axioms give a characterization
for the other  categories of bundles.
\end{remark}

We next show that $w_1(\xi) = 0$ if and only if $\xi$ is orientable. 

The oriented MacPhersonian $\Oo\MacP(k,n)$ is defined in~\cite{homotopy}.
The elements of the poset $\Oo\MacP(k,n)$
are all chirotopes of elements of $\MacP(k,n)$.
In~\cite{homotopy} it is shown that $\|\Oo\MacP(k,n)\|$ is
the universal double cover of $\|\MacP(k,n)\|$.  One can also define 
$\Oo\MacP(k,\infty)$ 
and show that its geometric realization is the double cover of
$\|\MacP(k,\infty)\|$.

\begin{defn}  An \textbf{orientation} of a matroid bundle $\xi = (B,\Mm)$
is a poset lifting 
$$\begin{array}{ccc}
 & &\Oo\MacP(k,\infty)\\
 &\nearrow&\downarrow\\
B&\stackrel{\Mm}{\rightarrow}&\MacP(k,\infty).
\end{array}$$
 \end{defn}

\begin{prop} Any topological lifting of $\|\Mm\|:\|B\|{\rightarrow}
\|\MacP(k,\infty)\|$ to $\|\Oo\MacP(k,\infty)\|$ is the geometric realization of
an orientation of $\Mm:B {\rightarrow} \MacP(k,\infty)$.
\end{prop}

\begin{proof} Any topological lifting is simplicial, and any simplicial lifting
to a poset covering space is the realization of a poset lifting. (This is
clear from looking at the lifting on individual simplices.)
\end{proof}

Thus an orientation of a matroid bundle is equivalent to an orientation of
the geometric realization of the associated combinatorial sphere bundle.
 
\begin{theorem} \label{orient}
Let $\xi$ be a vector bundle, a matroid bundle, or a
spherical  (quasi)-fibration. Then $\xi$ is orientable if and only if $w_1(\xi)=0$.
\end{theorem}

\begin{proof} Suppose first that $\xi=(B,\Mm)$ is a matroid bundle.  From the 
double cover  result, $H^1(\|\MacP(k,\infty)\|) \cong \Z2$, and is generated by
$w_1(\gamma_k)$ (where $\gamma_k$ is the universal bundle) since the first
Stiefel-Whitney class is a non-trivial characteristic class.  

Note that the map $\|\Oo\MacP(k,\infty)\|\rightarrow\|\MacP(k,\infty)\|$
is an $S^0$ bundle, and hence has a classifying map into $\RPinfty$. Thus we
have maps
$$\begin{array}{ccccc}
 &           &\|\Oo\MacP(k,\infty)\|&\rightarrow&             S^\infty\\
 &           &\downarrow   &                         &\downarrow\\
\|B\|&\stackrel{\|\Mm\|}\rightarrow&\|\MacP(k,\infty)\|   &
\stackrel{c}{\rightarrow}&\RPinfty
\end{array}
$$
Let $\beta:\|B\|\rightarrow\RPinfty$ be the composition of the lower two maps, 
 $\omega$ be the generator of $H^1(\RPinfty;\Z2)$.  Then by covering space theory 
$\beta$ has a lifting
if and only if $\beta^*\omega = 0$.  One can see directly that the
combinatorial vector bundle corresponding to the M\"obius strip mapping to the 
circle is non-orientable, and thus $c^*\omega\neq 0$.  Hence $c^*\omega =
w_1(\gamma_k)$, and the result follows. 

In the other cases, it suffices to consider a spherical fibration.  One could 
either use the classifying spaces $BSG(k)$ and $BG(k)$ for (oriented) spherical
fibrations and proceed as above, or use the fact the $Sq^1$ is the mod 2 Bockstein
to show that $w_1(\xi) =0$ if and only if the coefficients in the spectral
sequence used in the Thom isomorphism theorem are untwisted. 

\end{proof} 

Finally, we note that proof of the Whitney sum formula for matroid bundles also 
shows:

\begin{prop}  Let $\xi_1$ and $\xi_2$ be matroid bundles with orientations.  Then
$$e(\xi_1 \oplus \xi_2) = e(\xi_1)\cup e(\xi_2).
$$
\end{prop}

In particular, the Euler class is an unstable characteristic class.  Indeed,  if
$\epsilon$ is a trivial ($\Mm$ is constant) bundle of rank greater than zero, then
$e(\xi \oplus \epsilon) = e(\xi) \cup 0 = 0$.

\section{Vector bundles vs. matroid bundles: the Comparison Theorem}
\label{vb and cvb}

\begin{CT}  Let $B$ be a regular cell complex.  The composite of the natural 
transformations
$$V_k(B)\xrightarrow{C} M_k(B) \xrightarrow{\|E_0\|} Q_k(B)
$$
coincides
with the forgetful map given by deleting the zero section of a vector bundle.

\end{CT}

Thus the Stiefel-Whitney classes of the combinatorialization of a real  
vector bundle coincide with those of the original bundle. In particular, as a 
corollary we have Theorem~A. In addition, since for every realized rank $n$
oriented matroid $M$ the map
$G(k,\Rn)\to\|\MacP(k,\infty)\|$ factors as
$$G(k,\Rn)\to\|\Gamma(k, M)\|\to\|\MacP(k,\infty)\|,$$
and since $G(k,\Rn) \to G(k,\Rinfty)$ gives a split surjection on mod 2 cohomology, 
we have Theorem~B.

\begin{remark}  The Comparison Theorem could also be stated universally by saying that there are
maps
$$
BO(k) \to \|\MacP(k,\infty)\| \to BG(k)
$$
covered by maps of spherical quasifibrations on the universal sphere bundles.
\end{remark}

Let $M$ be a rank $n$ oriented matroid realized by a collection  $\{\phi_1,
\ldots,\phi_m\}$ of linear forms on $\Rn$.   Let
$S(k,\Rn)$ be the sphere bundle of the canonical bundle over
$G(k,\Rn)$.  An element of $S(k,\Rn)$ is a pair $(V,p)$ where $V$
is a $k$-plane in $\Rn$ and $p \in V$ has unit length.  Let $E_0(k,M)$ be the combinatorial sphere
bundle of the canonical bundle over 
$\Gamma(k,M)$.  Define  $\nu(V,p) =
(\mu(V),X(p)) \in E_0(k,M)$, where $\mu:G(k,\Rn)\to\Gamma(k,M)$ is the 
function defined in Section~\ref{defn:combgrass},
$X(p) = (\text{sign }\phi_1(p), \dots,
\text{sign }\phi_m(p))$ is a sign vector.
  Thus we have a commutative diagram 

$$
\begin{CD}
S(k,\R^n) @>\nu>> E_0(k,M) \\
@VV p V @VV\pi V \\
G(k,\Rn) @>\mu>> \Gg(k,M)
\end{CD}
$$
 with  $\mu$ and $\nu$ upper semi-continuous.    

\begin{lemma}  \label{lemma:comp}
Let $M$ be a realized rank $n$ oriented matroid.  Then there is a 
homotopy commutative diagram of continuous maps
$$
\begin{CD}
S(k,\R^n) @>\tilde \nu>>\| E_0(k,M)\| \\
@VV p V @VV\|\pi\| V \\
G(k,\Rn) @>\tilde \mu>> \|\Gg(k,M)\|
\end{CD}
$$
so that there is a $V\in G(k,\Rn)$ so that $\tilde \nu$ gives a homotopy 
equivalence from  the fiber above
$V$ to the fiber above
$\tilde \mu(V)$.  Furthermore $\tilde \mu$ is the map specified 
by Lemma \ref{Mac}.
\end{lemma}

\begin{proof} We again use Appendix~\ref{PL} and the semi-algebraic 
triangulation theorem of \cite{Hiro}. As in the proof of Lemma~\ref{Mac}, 
there is a semi-algebraic triangulation $T_G: \|L\|\to G(k,\Rn)$ refining the 
stratification of $G(k,\Rn)$ and a map $\tilde \mu:G(k,\Rn)\to\|\Gamma(k,M)\|$
so that $\tilde \mu\circ T_G$ is simplicial. Now consider the semi-algebraic
stratification of $S(k,\Rn)$ given by the intersections of
the preimages of elements under $\nu$ and the
preimages of simplices under $\Delta T_G^{-1}\circ p$. By \cite{Hiro}, there is a 
triangulation $T_S:\|K\|\to S(k,\Rn)$ refining this stratification, and so by
Corollary~\ref{lemma:triang} there is a map $\tilde\nu:S(k,\Rn)\to \|E_0(k,M)\|$
so that $\tilde \nu\circ \Delta T_S$ is simplicial.

To see that these maps make the above diagram commute up to homotopy,
let $s\in S(k,\Rn)$. Then $s$ lies in a simplex $\kappa \subset S(k,\Rn)$ of $
\Delta T_S$ 
and there is a simplex $\lambda$ of $\Delta T_G$ so that $p(\kappa) \subseteq 
\lambda \subset G(k,\Rn)$.
By construction $\tilde \nu (\kappa )$ is contained in the simplex spanned by the totally ordered 
set $\nu(\kappa)$ and $\tilde \mu (\lambda )$ is contained in the simplex spanned by the totally ordered 
set $\mu(\lambda)$.  Then $\tilde\mu(p(s))$ and $\|\pi\|(\tilde\nu(s))$ both lie 
in the closed simplex spanned by $\mu(\lambda)$, so there
is a straight-line homotopy from
$\|\pi\|\circ\tilde\nu$ to $\tilde\mu\circ p$.

Finally, note that for every vertex $V \in G(k,\Rn)$ of $ T_G$, the topological realization
theorem gives a homotopy equivalence from
the fiber over $V$ to the fiber over
$\tilde \mu(V)= \mu(V)$.
\end{proof}

\begin{lemma}  There is a homotopy
commutative diagram 
$$
\begin{CD}
S(k,\Rinfty) @>\tilde \nu>>\| E_0(k,\infty)\| \\
@VV p V @VV\|\pi\| V \\
G(k,\Rinfty) @>\tilde \mu>> \|\MacP(k,\infty)\|
\end{CD}
$$
so that there is a $V\in G(k,\Rinfty)$ so that $\tilde \nu$ gives a homotopy equivalence from  the
fiber above
$V$ to the fiber above
$\tilde \mu(V)$.  Furthermore $\tilde \mu$ is in the homotopy class of maps 
specified by Theorem~\ref{mu}.
\end{lemma}

\begin{proof}  The proof uses Theorem \ref{thm:triang} and the techniques of the proof of the 
previous lemma.
\end{proof}

\begin{proof}[Proof of the Comparison Theorem]  Let $\xi = (p:E \to B)$ be a vector bundle.  
Convert the map $\|\pi\|$ to a fibration $\pi_{\|\pi\|}$ and consider the diagram
$$
\begin{CD}
E_0(\xi) @>>> S(k,\Rinfty) @>\tilde \nu>>\| E_0(k,\infty)\| @>h>> P_{\|\pi\|}\\
@VVV @VV p V @VV\|\pi\| V @VV\pi_{\|\pi\|}V \\
B @>c>>G(k,\Rinfty) @>\tilde \mu>> \|\MacP(k,\infty)\| @>\text{Id}>> \|\MacP(k,\infty)\|
\end{CD}
$$

 By the homotopy lifting property, there is a map 
$\tilde
\nu'
\simeq h \circ
\tilde \nu$ so that $\pi_{\|\pi\|} \circ \tilde
\nu' = \tilde \mu \circ p$.  Furthermore, $h$ gives  a homotopy
equivalence on fibers (since $\|\pi\|$ is a quasifibration) and
$\tilde \nu$ gives a homotopy equivalence on a fiber,   so $\tilde \nu'$ gives a
homotopy equivalence on fibers.

Recall  $C[\xi]$ is defined by applying the
Simplicial Approximation Theorem  to find a subdivision  $B'$ of $B$ and a
map $c' : \Ff (B') \to \MacP(k,\Rinfty)$, where $\|c'\|$ is homotopic to  $\tilde
\mu \circ c$.  Then  $C[\xi] = [c']$.

We then have the following
equations in $Q_k(B)$:
\begin{align*}  \|E_0\| \circ C[\xi] &= \|E_0\|[c']\\
&= [\|E_0\| {c'}^*(\gamma^k)]\\
&= [\|{c'}^*E_0(\gamma^k)\|]\\
&=[{\|c'\|}^*\|E_0(\gamma^k)\|]\\
&= [{c^*}{\tilde \mu}^*\|E_0(\gamma^k)\|]\\
&=[{c^*}{\tilde \mu}^*P_{\|\pi\|}]\\
&= [c^*S(k,\Rinfty)]\\
&= [E_0(\xi)]\\
\end{align*}
\end{proof}

\section{Homotopy groups of the combinatorial Grassmannian} \label{sec:homotopy}

One can use the classical $J$-homomorphism to obtain limited information about
$\pi_i\|\MacP(k,n)\|$, or more generally about the homotopy groups of $\|\Gamma (k,M^n)\|$ for realizable $M^n$.    The idea
is use vector bundles over spheres to construct elements and use homotopy groups of spheres and characteristic classes to detect
them.  

A few remarks will give a context for these results.
The duality theorem for oriented matroids
gives $\|\MacP(k,n)\| \cong \|\MacP(n-k,n)\|$.  It is easy to show 
$\|\MacP(k,n)\|$ is connected, but there exist examples of $M^n$ such that
$\|\Gamma(n-1, M^n)\|$ is disconnected (\cite{MR93}).  
In \cite{homotopy}, it
was shown that $\pi_1  \|\MacP(k,n)\| \cong \pi_1G(k,\Rn)$, and stability 
results for large $n$ were established. There are also some results
on homotopy type of combinatorial Grassmannians for small values of $k$
or $n$
and for oriented matroids with few elements
(cf.~\cite{MZ}, \cite{Ba}, \cite{SZ}).

See~\cite{Whitehead} for the definition of the $J$-homomorphism 
$J_{i,k}: \pi_iO(k) \to \pi_{i+k}S^k$.
The limit as $k \to \infty$ is denoted
$$
J_i: \pi_iO \to \pi_i^S. 
$$
Stability results for the domain and range of $J$ show
\begin{align*}
\Im J_{i,k} \to & \Im J_{i,k+1}  &\text{is an epimorphism if $k \geq i+1$ and}\\
\Im J_{i,k} \to & \Im J_{i}  &\text{is an isomorphism if $k > i+1$}.
\end{align*}

A group $H$ is a {\bf subquotient} of a group $G$ if $H$ is isomorphic to a subgroup of a
quotient group of $G$.

\begin{theorem} \label{thm:j}
 Let $M^n$ be a realized
rank $n$ oriented 
matroid. Let $p$ be a point in the image of $\tilde\mu: G(k,\Rn)\to\|\Gamma(k,M^n)\|$.
\begin{enumerate}
\item \text{$\Im J_{i-1,k}$  is a subquotient of $\pi_i(\|\Gamma(k,M^n)\|,p)$  when $n-k \geq i$.}
\item \text{$\Im J_{i-1}$  is a subquotient of $\pi_i(\|\Gamma(k,M^n)\|,p)$  when $n-k \geq i$ and $k>i$.}
\end{enumerate}
\end{theorem}

\begin{proof}  Let $G(k)$ denote the monoid of self-homotopy equivalences of
$S^{k-1}$, given the compact-open topology.  Its classifying space $BG(k)$ classifies rank $k$
spherical (quasi)-fibrations \cite{S}. $G_0(k+1)$ denotes the monoid of self-homotopy equivalences
of $S^{k}$ which fix a base point $*$. 

The result follows from commutativity up to homotopy
of the  diagram
$$
\begin{CD}  G(k,\Rn) @>\tilde\mu>> \|\Gamma(k,M^n)\| \\
@VV\beta V @VV\delta V \\
G(k,\Rinfty) @>\gamma>> BG(k) @>\epsilon >> BG_0(k+1),
\end{CD}
$$
the surjectivity of $\pi_i (\beta)$
when $n-k\geq i$, and the identification of
$\pi_i(\epsilon
\circ \gamma)$  with $J_{i-1,k}$.  
The map $\beta$ is given by inclusion, the map $\gamma$ by classifying the canonical bundle minus
its zero section, and the map $\epsilon$ is $B$ applied to the injection of monoids $G(k) \to
G_0(k)$ given by suspension.  The map $\tilde\mu$, the map $\delta$, and the homotopy 
$\delta \circ \tilde\mu \simeq \gamma \circ \beta$, are given by our three main theorems:
the Combinatorialization Theorem, Spherical Quasifibration Theorem, and the Comparison Theorem (see also Lemma
\ref{lemma:comp}).  

The surjectivity of $\pi_i(\beta)$ when $n-k \geq i$ follows either from the Cellular Approximation
Theorem applied to the Schubert cell decomposition, or by viewing $G(k,n)$ as $O(n)/(O(k) \times
O(n-k))$ and using the homotopy long exact sequence of a fibration.  

Now $\pi_iG(k,\Rinfty)\cong \pi_{i-1}O(k)$ since there is a fibration $O(k) \to V(k,\Rinfty) \to
G(k,\Rinfty)$ and the Stiefel manifold is contractible. Also $\pi_iBG_0(k+1) \cong \pi_i G_0(k+1)$,
which is in turn $\pi_{i+k}S^k$ by the adjoint property of smash and mapping spaces in the category
of based CW complexes.  Thus we have identified the domain and range of $\pi_i(\epsilon
\circ \gamma)$ with that of $J_{i-1,k}$, and the identification of the two maps consists of tracing
through these identifications.
\end{proof}

\begin{corollary}  Let $M^n$ be 
a realized rank $n$ oriented 
matroid. Let $p$ be a point in the image of $\tilde\mu: G(k,\Rn)\to\|\Gamma(k,M^n)\|$.
\begin{enumerate} 
\item  
$$\pi_2(\|\MacP(k,n)\|) \cong \pi_2(G(k,\R^n))  \cong \begin{cases}
0 & \text{if $k=1$}\\
{\mathbb Z} & \text{if $k = 2$ and $n-k \geq 2$}\\
\Z2  & \text{if $k \geq 3$ and $n-k \geq 3$.}
\end{cases} $$
\item  $\pi_k (\|\Gamma(k,M^n)\|,p)$ has $\Zz$ as a subgroup when $k$ is even and $n \geq 2k$.
\item  $\pi_i (\|\Gamma(k,M^n)\|,p)$ has $\Z2$ as a subquotient when $i \equiv 1,2$  (mod 8), $n-k \geq i$, and $k \geq i$. 
\item  $\pi_{4m} (\|\Gamma(k,M^n)\|,p)$ has $\Zz_{a_m}$ as a subquotient when $m > 0$, $n-k \geq 4m$, and $k \geq 4m$.  Here $a_m$
is the denominator of $B_m/4m$ expressed as a fraction in lowest terms, and $B_m$ is the $m$-th Bernoulli number.
\end{enumerate}
\end{corollary}

\begin{proof}  1.  By \cite{homotopy}, the combinatorialization map $\pi_2(G(k,\R^n)) \to \pi_2(\|\MacP(k,n)\|)$ is surjective, so
(1) follows from (2) and (3).  An alternate route to (1) is to use characteristic classes to detect elements of
$\pi_2\|\MacP(k,n)\|$, by applying the Euler class and second Stiefel-Whitney class of the combinatorialization of the complex
Hopf bundle over the 2-sphere.  
 
2. There is a characteristic class version and a homotopy theoretic 
version of the proof.
  The characteristic class proof is as follows. Consider
the oriented combinatorial Grassmannian $\tilde\Gamma(k,M^n)$, defined analogously to $\Oo\MacP(k,n)$. The forgetful map
$\|\tilde\Gamma(k,M^n)\|\to
\|\Gamma(k,M^n)\|$ is a double cover, so $\pi_k(\|\Gamma(k,M^n)\|,p)
\cong
\pi_k(\|\tilde\Gamma(k,M^n)\|,p')$ for each lifting $p'$ of $p$.
The tangent bundle of the $k$-sphere combinatorializes to a map $S^k\to \tilde G(k,\Rn) \xrightarrow{\tilde \mu} 
\|\tilde\Gamma(k,M^n)\|$. The evaluation of the Euler class of the 
tangent bundle of a manifold on the fundamental class of that manifold is the Euler characteristic of the manifold (cf. 11.12 in~\cite{MS}). 
In particular, the Euler class of the
tangent bundle of an even-dimensional sphere represents twice the generator of 
$H^k(S^k)\cong\Zz$. The Euler class applied to oriented vector bundles 
over $k$-spheres can be viewed as a homomorphism from 
$\pi_k(BSG(k))$ to $H^k(S^k)$, where $BSG(k)$ classifies  rank $k$ oriented spherical fibrations. Thus the combinatorialization of
this tangent bundle
generates an infinite subgroup of $\pi_k(\|\tilde\Gamma(k,M^n)\|,p'))$.

The homotopy theoretic version is to consider the Hopf invariant
$$
H : \Im J_{k-1,k} \to \Zz.
$$
The domain of $ J_{k-1,k}$ is $\pi_{k-1}(O(k))$, which classifies 
$k$-bundles over $S^k$.
Now $H \circ J_{k-1,k}$ is the Euler class, so $2 \Zz \subseteq \Im H$, 
using the tangent bundle of the $k$-sphere again.  (For the reader's edification, we note the Hopf
invariant is onto if and only if $k = 2, 4$ or 8, as can be seen by Bott periodicity or by Adams's work on Hopf invariant one.)

3.,4.  This follows from Theorem~\ref{thm:j} and
the deep homotopy theoretic computation of $\Im J$ due to Adams \cite{adams}.  The result is that $\Im
J_{i-1}$ is $\Z2$ for $i \equiv 1,2$ (mod 8), $\Zz_{a_m}$ for $i = 4m$, and is zero otherwise.
\end{proof}

An exposition of Bernoulli numbers and topology is given in \cite[Appendix B]{MS}.  The first few values of $a_m$ are :

\begin{tabular}{c|c|c|c|c|c|c|c}
$a_1$&$a_2$&$a_3$&$a_4$&$a_5$&$a_6$&$a_7$&$a_8$\\ \hline
24& 240&504&480&264&65520&24&16320
\end{tabular}


By the stability results of~\cite{homotopy}, the results
of our corollary also apply to the MacPhersonian $\MacP(k,\infty)$.

\section{Vector fields and characteristic classes}\label{obstr}

The classical motivation for characteristic classes  was as obstructions to the
existence of linearly independent vector fields of a manifold, or more
generally, independent sections of a vector bundle.  This section gives 
 a combinatorial analog.

\begin{defn}  The \textbf{combinatorial Stiefel manifold}  $V_l(k,n)$
is the poset of all $M\in\MacP(k,n+l)$ satisfying both of the conditions below
\begin{enumerate}
\item $\text{rank} (M \backslash \{n+1,\ldots,n+l\}) = k$
\item $\{n+1,\ldots,n+l\}$ is independent in $M$.
\end{enumerate}
\end{defn}

Note that deleting $\{n+1,\ldots,n+l\}$ gives a surjective map
$V_l(k,n)\rightarrow\MacP(k,n)$ if $l \leq k$.

\begin{defn}
If $(B,\xi)$ is a matroid bundle, an \textbf{independent set of $l$ 
vector fields} is a lifting
$$\begin{array}{ccc}
 & &V_l(k,n)\\
 &\nearrow&\downarrow\\
B&\stackrel{\xi}{\rightarrow}&\MacP(k,n).
\end{array}$$
 \end{defn}

\begin{lemma} If a matroid bundle $\xi=(B, \Mm: B\to\MacP(k,n)$ 
admits a set $\nu:B\rightarrow V_l(k,n)$ of $l$ independent vector fields,
then:
\begin{enumerate}
\item The map 
$$\begin{array}{rrcl}
\Qq:&B&\rightarrow&\MacP(k-l,n)\\
    &\sigma&\mapsto&\nu(\sigma)/\{n+1,\ldots,n+l\}
\end{array}$$
is a matroid bundle.

\item If $\epsilon_l$ is the trivial rank $l$ bundle over $B$
sending each cell to the rank $l$ oriented matroid with elements 
$\{n+1,\ldots,n+l\}$, then the matroid bundles $\xi$ and  $\Qq \oplus\ \epsilon_l$ are
$\|B\|$-isomorphic.
\end{enumerate}
\end{lemma}

$\Qq$ is called the \textbf{quotient bundle} of $\nu$.

\begin{proof} (1) This follows immediately from Lemma~\ref{contract}

(2) The proof is by induction on $l$, and relies on the Order Homotopy
Lemma.

If $l=1$, first note that $\nu(\sigma)\geq\xi(\sigma)$ for all $\sigma$,
so $\nu\simeq \xi$.  (Here $\simeq$ means $\|B\|$-isomorphic.   All homotopies
occur  in
$\|\MacP(k,n+l)\|$; there is no need to go to $\infty$.)  Thus,  it suffices to
show that
$\nu(\sigma)\geq
\Qq(\sigma)\oplus\epsilon_1(\sigma)$ for all $\sigma$.  Any covector  $X \times Y
\in \Vv^* (\Qq(\sigma)\oplus\epsilon_1(\sigma))$ is built from covectors $X \in
\Vv^*(\Qq(\sigma)) \subset
\Vv^*(\nu(\sigma))$ and
$Y:
\{n+1\}
\to
\{+,-,0\}$. Since $\{n+1\}$ is independent in $\nu(\sigma)$, there is a 
$\widetilde{Y} \in
\Vv^*(\nu(\sigma))$, so that $\widetilde{Y} \geq Y$.  Thus $X \circ \widetilde{Y}$ 
is a covector of
$\nu(\sigma)$ which is greater than or equal to $X \times Y$.

For  $l > 1$, let $\nu'$ be the set of $l-1$ independent vector fields 
obtained from $\nu$ by deleting $n+l$ from each oriented matroid $\nu(\sigma)$,
and  let $\Qq'$ be the resulting quotient bundle. Note that the vector field
$\nu_l$ on $\Qq'$ given by $n+l$ is non-vanishing and $\Qq$ is the quotient 
bundle of $\Qq'$ by
$\nu_l$.  Hence we have
\begin{align*} \xi &\simeq \Qq' \oplus \epsilon_{l-1} &\text{by the induction 
hypothesis} \\ &\simeq \Qq \oplus \epsilon_1 \oplus \epsilon_{l-1} &\text{by the
$l=1$ case.}
\end{align*}

\end{proof}

\begin{theorem}\label{vectorfields}
 If a rank $k$ matroid bundle $\xi$ admits an independent set of $l$
vector fields, then $w_{k-l+1}(\xi)=0$.
\end{theorem}

\begin{proof}[Proof of Theorem~\ref{vectorfields}]
 If $\nu$ is a set of $l$ independent vector fields in $\xi$ and $\Qq$ is the 
resulting quotient bundle, then by the above lemma,
$\xi\simeq\Qq\oplus\epsilon_l$. Thus by the  Whitney sum formula,
$w(\xi)=w(\Qq)w(\epsilon_l)$ where $w = 1 + w_1 +w_2 + \cdots$ is the total Stiefel-Whitney class.
But
$w(\epsilon_l)=1$ since
$\epsilon_l$ is trivial, and $w_{k-l+1}(\Qq)=0$ since $\Qq$ is a rank $k-l$ bundle.
\end{proof}

Similarly, we have

\begin{theorem} If a rank $k$ matroid bundle $\xi$ admits a non-zero cross section 
(i.e. an independent set of 1 vector field) then the Euler class $e(\xi)$ vanishes.
\end{theorem}

\section{Some open questions}

There are open questions everywhere you spit; we list a few.

\begin{enumerate}
\item  Is a CD manifold a Poincar\'e complex? Does a CD manifold satisfy
Poincar\'e duality?

\item  Give a definition of isomorphism of CD-manifolds; show that a diffeomorphism class of smooth
manifolds determines an isomorphism class of CD-manifolds.

\noindent {\em Discussion:  Problems (1) and (2) are manifold theoretic analogues of the bundle
theoretic results of this paper.  Macpherson \cite{Mac} defined a CD manifold
and showed how a smooth manifold with a smooth triangulation determines a
CD manifold.  He asked whether a CD-manifold was a topological manifold; this
was shown in a special case in \cite{CDmanifolds}.}

\item Are there exotic mod 2 characteristic classes?

\item Are there Pontrjagin classes?

\noindent {\em Discussion:  The authors together with Eric Babson have outlined a construction of
rational Pontrjagin classes.  If these classes were integral cohomology classes, that would imply
the existence of exotic CD 7-spheres.  
}

\item  Is $\MacP(\infty,\infty)$ an infinite loop space?

\item  Compute the homotopy groups of $\MacP(\infty,\infty)$.

\noindent {\em Discussion:  Solving  questions 5 and 6 would show that combinatorial vector bundles
give a generalized cohomology theory and compute the coefficient groups (isomorphism classes of
matroid bundles over spheres.)}
\end{enumerate}

\appendix

\section{Topological maps from combinatorial ones} \label{PL}

Section~\ref{defn:combgrass} described a natural map 
$\mu : G(k,\Rn) \to \Gamma(k,M)$ for any rank $n$ oriented matroid $M$ with a
fixed realization in $\Rn$, given by intersecting the hyperplanes of the
realization with
$V \in G(k,\Rn)$ and taking the corresponding oriented matroid. 
This appendix describes how we use this map to make a simplicial map from a
triangulation of $G(k,\Rn)$ to $\Delta \Gamma(k,M)$, unique up to homotopy.
We also construct a topological map $G(k,\Rinfty)
\to\|
\MacP(k,\infty)\|$ which is in some sense the limit of the simplicial maps obtained when 
$\Gamma(k,M) = \MacP(k,n)$. As discussed in~\cite{notPL}, this limit map is not PL.

 Actually, we will
work more generally, considering  maps from spaces to posets satisfying certain
properties.

\begin{defn}  A \textbf{triangulation} of a topological space $X$ is a
homeomorphism $T : \|K\| \to X$ where
$K$ is a simplicial complex.  We will abuse language slightly and refer to the image 
$T(\|\sigma\|)$ of a simplex $\|\sigma\|$ under a triangulation $T$ as
a \textbf{simplex of $T$}. If $X$ resp. $Y$ are spaces equipped with
triangulations $S$ resp. $T$
then a map $f : X \to Y$ is {\bf simplicial} if $T^{-1} \circ f \circ S$ is.  A
\textbf{subdivision} of a simplicial complex
$K$ is a homeomorphism $S : \|K'\|
\to  \|K\|$ where $K'$ is a simplicial complex and for every simplex $\sigma' \in K'$, there is a simplex
$\sigma \in K$ so that $S(\|\sigma'\|) \subseteq \|\sigma\|$ and $S$ is linear on 
$\|\sigma'\|$.  The triangulation $T'  = T \circ S$
 is a \textbf{subdivision of the 
triangulation
$T$}.   An example of such  is the barycentric subdivision $\Delta T : \| \Delta K
\| \to X$.  Two triangulations of a space
$X$ are
\textbf{equivalent} if they have a common subdivision.  A \textbf{PL space} is a
space equipped with a fixed equivalence class of triangulations, called
the \textbf{PL triangulations}. A map $f : X\to
Y$  between PL spaces is a \textbf{PL map} if there are PL triangulations 
so that $f$ is simplicial.
\end{defn}

\begin{defn}  \label{strat} Let $\mu : G \to M$ be a function from a space to a poset. 
The partition $\{\mu^{-1}(m):m\in M\}$ of $G$ is the \textbf{stratification of $G$
induced by $\mu$}. The map $\mu$ is
\textbf{upper semi-continuous} if every $g \in G$ has a neighborhood $U$ so that $\mu(U) \subseteq
M_{\geq \mu(g)}$. Thus the closure of a stratum $\mu^{-1}(m)$ maps to the 
lower order ideal $M_{\leq m}$ of $M$. A triangulation
$T : \|K\| \to G$ \textbf{refines the stratification} if the interior of every 
simplex maps under  to a single element of
the poset $M$.  
\end{defn}

\begin{lemma} \label{ordered}
 Let $\mu : G \to M$ be upper semi-continuous and let $T : \|K\| \to G$ be
a triangulation refining
the stratification.  Then for any simplex $\alpha$ of the barycentric
subdivision,
$\mu(\Delta T \|\alpha \|)$ is totally ordered.
\end{lemma}

\begin{proof}  For a simplex $\sigma$ of $K$, let $\langle \sigma \rangle$ 
denote both the barycenter of
the simplex $\|\sigma \| \subseteq \|K\|$ and the corresponding vertex of $\|\Delta K\|$.  Let
$\alpha = \{\sigma_n >
\dots > \sigma_0\} \in \Delta K$ be a chain of simplices of $K$.  Note that 
$$\|\alpha \| \subseteq (\text{int } \|\sigma_n\|) \cup \|\{\sigma_{n-1} >
\dots > \sigma_0\}\|
$$
hence inductively, 
$$
\mu(\Delta T \|\alpha \|) =
\{\mu(T\langle \sigma_n\rangle), \dots,\mu(T\langle \sigma_0\rangle)\}.
$$ 
Since the triangulation
refines the stratification and $\mu$ is upper semi-continuous, 
$\mu(T\langle \sigma_i\rangle) \geq \mu(T\langle \sigma_{i-1}\rangle)$ for all $i$.    
\end{proof}

\begin{corollary}  \label{lemma:triang} Let $\mu : G \to M$ be upper semi-continuous and let $T : \|K\| \to G$ be
a triangulation refining the stratification.
\begin{enumerate}  
\item  There is a map $\tilde \mu_T : G \to \|M\|$, simplicial with respect to the barycentric 
subdivision of $T$, which agrees with $\mu$ on the vertices of $\Delta T$.  
\item  If $T'$ is a subdivision of $T$, then $\tilde\mu_T \simeq\tilde\mu_{T'}$.
\end{enumerate}
\end{corollary}

\begin{proof}  Part (1) follows from the last lemma by defining $\tilde\mu_T$ on vertices 
and extending by linearity with respect to $\|K\|$.  Part (2) follows from a straight-line homotopy 
$t\tilde \mu_T(a) + (1-t) \tilde \mu_{T'}(a)$ across the simplices of $\|M\|$.  
\end{proof}

We need an infinite version of the corollary.

\begin{theorem}\label{thm:triang}
Let $G_1\subseteq G_2\subseteq\cdots$ be a sequence of PL spaces, $M_1\subseteq M_2\subseteq\cdots$
a sequence of posets, and $\mu_1:G_1 \to M_1, \mu_2:G_2 \to M_2, \dots $ a sequence of upper
semi-continuous maps so that $\mu_{i}|_{G_{i-1}}=\mu_{i-1}$.  Let $T_1, T_2, \dots$ be a sequence of PL
triangulations of $G_1, G_2, \dots$ refining the stratifications given by $\mu_1, \mu_2, \dots$, 
so that the restriction of the $i$-th triangulation is a subdivision of the $(i-1)$-st triangulation.  
Let $M$ be the union of the $M_i$'s and $G$ be the direct limit of the $G_i$'s.\
\begin{enumerate}
\item There is a continuous map $\tilde\mu : G \to \|M\|$ which, for every $i$,
restricts to a map $\tilde\mu_i:G_i\to\|M_i\|$ which is simplicial with respect to the 
barycentric subdivision of $T_i$, and so that $\tilde\mu_i$ agrees with $\mu_i$ on the vertices of 
$\Delta T_i$ in $G_i \backslash G_{i-1}$.
\item Let $T'_1, T'_2, \dots$ be another sequence of PL triangulations 
satisfying the same hypotheses as $T_1, T_2, \dots$. Let 
$\tilde\mu': G \to \|M\|$ be the map given by part (1) using the $T'_i$'s.
  Then $\tilde\mu \simeq \tilde\mu'$ and the 
homotopy restricts to a homotopy $\tilde\mu_i \simeq \tilde{\mu'}_i$ for all $i$.  
\end{enumerate}
\end{theorem}

\begin{proof}  We will inductively define maps $S_i : G_i \to \Delta M_i$ so that 
\begin{enumerate}
\item  $S_i|_{G_{i-1}} = S_{i-1},$
\item For a vertex $v \in G_i \backslash G_{i-1}$ of $\Delta T_i$,
$$S_i(v) = \{\mu(v)\},
$$
\item  For a point $p \in G_i$ in the interior of a simplex of $\Delta T_i$ which is spanned by 
vertices $\{v_0, \dots, v_n\} \subset G_i$,
$$S_i(p) = \cup_j S_i(v_j).
$$
\end{enumerate}

Assume inductively that $S_{i-1}$ has been defined satisfying (1), 
(2), and (3), and also assume 
inductively that $S_{i-1}$ satisfies property 

\vskip .1in

\noindent
4.  $\text{max }S_{i-1}(p)=\mu(p)  \qquad \text{for }p \in G_{i-1}$.
\vskip .1in

    We then use property (1) to define $S_i$ on $G_{i-1}$ and properties (2)
 and (3) to define
$S_i$ on $G_{i}\backslash G_{i-1}$.  We need to verify three things: first that property (3) 
continues to hold for $S_i$ and $p \in G_{i-1}$, second that $S_i(p)$
is totally ordered for $p \in G_{i}\backslash G_{i-1}$, and third that property (4) holds.
We leave the verification of the first part to the reader.

We now show that $S_i(p)$ is totally ordered for $p \in G_{i}\backslash
G_{i-1}$.  Suppose $p$
is in the interior of a simplex of $\Delta T_i$ with vertices
\[
\{v_0, \dots, v_n\} \cup \{w_0, \dots, w_m\}
\]
with the $v$'s in $G_{i-1}$ and the $w$'s in $G_{i} \backslash G_{i-1}$.
Then by
choosing a point $q$ in the interior of the simplex spanned by the $v$'s
and by the (omitted) proof
that property (3) holds for $S_i$ and $q$, we see $\cup_j S_i(v_j)$ is
totally ordered.
The proof of Lemma \ref{ordered} shows that $\cup_j S_i(w_j)=\{\mu(w_0),
\dots, \mu(w_m) \}$ is totally ordered.
Because $G_{i-1}$ is a subcomplex of the triangulation $T_i$, because we
have taken a barycentric
subdivision, because $\mu$ is upper semi-continuous, and because $T_i$
refines the
stratification, $\mu(w_j) \geq \mu(v_k)$ for all $j$ and $k$.  Hence    
$S_i(p)$ is totally ordered.
Finally property (4) holds since we have taken a barycentric subdivision.

Next we use $S_i$ to define $\tilde\mu_i$ inductively.  For $p \in G_{i-1}$,
let
$\tilde\mu_i(p) = \tilde\mu_{i-1}(p)$.  For a vertex
$p \in G_{i}\backslash G_{i-1}$ of $\Delta T_i$, define $\tilde\mu_i(p)=
\mu(p)$.  Then for $p$ in the
interior of a simplex spanned by $\{v_0, \dots, v_k\}$, define
$\tilde\mu_i(p)$ by linearity, noting
inductively that
$\mu_i(p)$ is in the closed simplex spanned by $S_i(p)$.  This completes
the proof of part (1) of
the Theorem.

For part (2), it suffices to consider the case when each $T_i'$ subdivides 
$T_i$.  Then note that $S'_i(p) \subseteq S_i(p)$, so we can use the
straight-line homotopy across the simplex spanned by $S_i(p)$.

\end{proof}

\section{Babson's criterion}  \label{eric}
This appendix gives the criterion of Babson for the realization of a 
poset map to
be a quasifibration.  We first state two results of Quillen \cite[page 98]{Q73}, specialized from
general categories to posets.  

\begin{QTA}  Let $f : P \to Q$ be a poset map.  If for all $q \in Q$, 
$\|f^{-1}(Q_{\geq q})\|$ is contractible,
then $\|f\|: \|P\| \to \|Q\|$ is a homotopy equivalence.
\end{QTA}

A commutative square of spaces
$$
\begin{CD}  A @>\alpha>> B \\
@V\beta VV @VV\gamma V \\
C @>\delta>> D
\end{CD}
$$
is {\bf homotopy cartesian} if for all $c \in C$, the induced map on homotopy fibers
$$\pi_{\beta}^{-1}(c) \to \pi_{\gamma}^{-1}(\delta c) 
$$
is a homotopy equivalence.

\begin{QTB} 
Let $f : P \to Q$ be a poset map such that for every inequality $q \geq q'$ in  $Q$, the inclusion
map
$\| f^{-1}(Q_{\geq q})\|\to \| f^{-1}(Q_{\geq q'})\|$  is a homotopy
equivalence. Then for any $q
\in Q$, 
$$
\begin{CD}
\|f^{-1}(Q_{\geq q})\| @>>> \|P\| \\
@VVV @VVV \\
\|Q_{\geq q}\| @>>> \|Q\|  
\end{CD}
$$
is homotopy cartesian.
\end{QTB}

Since the geometric realization of a poset can be identified with the realization 
of the poset obtained by reversing the inequalities, one may reverse
 the inequalities
in Quillen's Theorem A and B.  Similar remarks apply to the various
lemmas below.

In Babson's thesis, both of the previous two results of Quillen were recast.  
\begin{lemma}  \label{lemma A:babson} If $f : P \to Q$ is a poset map satisfying 
both of the conditions below, then
$\|f\|$ is a homotopy equivalence.
\begin{enumerate}
\item $\| f^{-1}q\|$ is contractible for all $q \in Q$.
\item $\|f^{-1}q \cap P_{\leq
p}\|$ is contractible whenever $p\in P$, $q \in Q$, and $q \leq f(p)$.
\end{enumerate}
\end{lemma}

\begin{proof} By Theorem A applied to $f$ and condition 1, we only need show that the
realization of  the inclusion $i : f^{-1}q
\to f^{-1}(Q_{
\geq q})$ is a homotopy equivalence for all $q \in Q$.  But this follows from 
condition 2 and by applying Theorem A 
to $i$,  noting that
$$i^{-1}(f^{-1}(Q_{
\geq q})_{\leq p}) = f^{-1}q \cap P_{\leq p}$$   
\end{proof}

Recall that a chain in a poset $P$ is a non-empty, 
totally ordered,  finite subset of $P$, and 
$\Delta P$ denotes the poset of chains in $P$, where the
partial order is given by inclusion.  Note 
$\|\Delta P\|$ is the barycentric subdivision of
$\|P\|$.  If $f : P \to Q$ is a poset map, then there is a poset map $\Delta f
: \Delta P \to \Delta Q$ sending a chain $c$ to the chain $f(d)$.  Note that
for a chain $d$ of $Q$, the symbols $f^{-1}d$ and $(\Delta f)^{-1}d$ have
different meanings.

\begin{lemma}\label{lemma B:babson}
If $f : P \to Q$ is a poset map satisfying the conditions below, then
$\|f\|$ is a quasifibration.
\begin{enumerate}
\item $\|f^{-1}q \cap P_{\leq
p}\|$ is contractible whenever $p\in P$, $q \in Q$, and $q \leq f(p)$.
\item $\|f^{-1}q \cap P_{\geq
p}\|$ is contractible whenever $p\in P$, $q \in Q$, and $q \geq f(p)$.
\end{enumerate}
\end{lemma}

\begin{proof}  We wish to  apply Quillen's Theorem B to the induced map
of posets
$\Delta f : \Delta P \to \Delta Q$.  For a chain $d \in \Delta Q$, 
there is a retraction
\begin{align*}
r : (\Delta f)^{-1} (Q_{\geq d}) & \to (\Delta f)^{-1} d \\
c &\mapsto c \cap f^{-1}d.
\end{align*}
By 
the order homotopy lemma,  $\|r\|$ is actually a deformation retraction  since 
$r(c)\leq c$ for every $c\in (\Delta f)^{-1}(Q_{\geq d})$.

For an inequality
$d
\geq d'$ in
$\Delta Q$, consider the following commutative diagram:
$$
\begin{CD}  (\Delta f)^{-1} d @>>>  (\Delta f)^{-1} (Q_{\geq d})\\
@V\alpha VV  @VVV  \\
(\Delta f)^{-1} d'   @<r<< (\Delta f)^{-1} (Q_{\geq d'})
\end{CD}
$$
where $\alpha (c) = c \cap f^{-1}d'$, $r$ is the retraction and the other two arrows are inclusions.  To verify the
hypothesis of Quillen's Theorem B applied to $\Delta f$, it suffices to prove the following lemma.

\begin{lemma}  Let $f : P \to Q$ be a poset map satisfying the conditions of Lemma \ref{lemma B:babson}.  Given an
inequality $d \geq d'$ in $\Delta Q$, the geometric realization of the map
\begin{align*}  \alpha(d \geq d') : (\Delta f)^{-1}d& \to (\Delta f)^{-1}d'\\
c &\mapsto c \cap f^{-1}d'
\end{align*}
is a homotopy equivalence.
\end{lemma}

\begin{proof} 
To prove this it suffices to check only those 
inequalities in $\Delta Q$ given by deleting the smallest or largest element
of a chain.
  Let $d = (q_r > \dots > q_1 >
q_0)$ be a chain in
$\Delta Q$.
Let $d' = (q_r > \dots > q_2 > q_1)$ and $d'' = (q_{r-1} > \dots > q_1 >
q_0)$.  We will use Lemma
\ref{lemma A:babson} to prove that $\|\alpha(d \geq d')\|$ and $\| \alpha(d \geq d'')\|$ are homotopy equivalences. 
Note
$\alpha(d \geq d')^{-1}c'$ is isomorphic to $\Delta (f^{-1}q_0 \cap
P_{\leq \text{min }c'})$ by the isomorphism $c \mapsto c \cap f^{-1}q_0$,
and the geometric realization of $\Delta (f^{-1}q_0 \cap
P_{\leq \text{min }c'})$ is contractible by condition (1).  

Now if $c' \subseteq c \cap
f^{-1}d'$ where $c$ and $c'$ are chains which map to $d$ and $d'$ respectively, 
$$\alpha(d \geq d')^{-1}c' \cap ((\Delta f)^{-1}d)_{\leq c}
$$
has a maximum element, namely $(c \cap f^{-1}q_0) \cup c'$, so its
geometric realization is
contractible. Thus by Lemma \ref{lemma A:babson}, $\|\alpha(d \geq d')\|$ is a homotopy equivalence.  The
proof that $\|F(d \geq d'')\|$ is a homotopy equivalence is similar and
uses property (3) of Lemma
\ref{lemma B:babson}.
\end{proof}
 
Now we return to the proof of Lemma \ref{lemma B:babson}.  Note that for $d \in \Delta Q$, the geometric realization
of the diagram
$$
\begin{CD} (\Delta f)^{-1}d @>>> (\Delta f)^{-1} (Q_{\geq d}) \\
@VVV @VVV \\
d @>>> Q_{\geq d} 
\end{CD}
$$ 
is homotopy cartesian.  Thus for all vertices and barycenters of $\|Q\|$ the inclusion of the fiber of $\|f\|:\|P\|
\to \|Q\|$ in the homotopy fiber is a homotopy equivalence.  Since this is true for the barycenters and since $\|f\|$
is a simplicial map, this is true for all points in the interior of a simplex, so $\|f\|$ is a quasifibration.
\end{proof}

\bibliographystyle{alpha}
\bibliography{biblio}

\vfill

\emph{\begin{tabular}{l}
Department of Mathematics\\
Texas A\&M University   \\
College Station, TX 77843\\[3pt]
laura.anderson@math.tamu.edu\\
http://www.math.tamu.edu/\mbox{$\sim$}laura.anderson/ 
       \end{tabular}}

\vskip 24pt
\emph{\begin{tabular}{l}
Department of Mathematics\\
Indiana University\\
Bloomington, IN 47405\\[3pt]
jfdavis@indiana.edu\\
http://www.indiana.edu/\mbox{$\sim$}jfdavis/ 
       \end{tabular}}

\end{document}

%% file: notbundle.pstex_t
\begin{picture}(0,0)%
\special{psfile=notbundle.pstex}%
\end{picture}%
\setlength{\unitlength}{3947sp}%
\begingroup\makeatletter\ifx\SetFigFont\undefined
\def\x#1#2#3#4#5#6#7\relax{\def\x{#1#2#3#4#5#6}}%
\expandafter\x\fmtname xxxxxx\relax \def\y{splain}%
\ifx\x\y   
\gdef\SetFigFont#1#2#3{%
  \ifnum #1<17\tiny\else \ifnum #1<20\small\else
  \ifnum #1<24\normalsize\else \ifnum #1<29\large\else
  \ifnum #1<34\Large\else \ifnum #1<41\LARGE\else
     \huge\fi\fi\fi\fi\fi\fi
  \csname #3\endcsname}%
\else
\gdef\SetFigFont#1#2#3{\begingroup
  \count@#1\relax \ifnum 25<\count@\count@25\fi
  \def\x{\endgroup\@setsize\SetFigFont{#2pt}}%
  \expandafter\x
    \csname \romannumeral\the\count@ pt\expandafter\endcsname
    \csname @\romannumeral\the\count@ pt\endcsname
  \csname #3\endcsname}%
\fi
\fi\endgroup
\begin{picture}(5355,1588)(976,-3875)
\put(3025,-2676){\makebox(0,0)[lb]{\smash{\SetFigFont{8}{9.6}{rm}$a$}}}
\put(2878,-3703){\makebox(0,0)[lb]{\smash{\SetFigFont{8}{9.6}{rm}$c$}}}
\put(3074,-3165){\makebox(0,0)[lb]{\smash{\SetFigFont{8}{9.6}{rm}$b$}}}
\put(6331,-2703){\makebox(0,0)[lb]{\smash{\SetFigFont{8}{9.6}{rm}$a$}}}
\put(6184,-3730){\makebox(0,0)[lb]{\smash{\SetFigFont{8}{9.6}{rm}$c$}}}
\put(5476,-3286){\makebox(0,0)[lb]{\smash{\SetFigFont{8}{9.6}{rm}$0=b$}}}
\put(976,-2761){\makebox(0,0)[lb]{\smash{\SetFigFont{8}{9.6}{rm}$M_1=$}}}
\put(4276,-2761){\makebox(0,0)[lb]{\smash{\SetFigFont{8}{9.6}{rm}$M_0=$}}}
\end{picture}

%% file: MacPherson.bbl
\newcommand{\etalchar}[1]{$^{#1}$}
\begin{thebibliography}{BLS{\etalchar{+}}93}

\bibitem[AD]{notPL}
L.~Anderson and J.~F. Davis.
\newblock {S}chubert stratifications and triangulations of the infinite
  {G}rassmannian.
\newblock To appear.

\bibitem[Ada65]{adams}
J.~F. Adams.
\newblock On the groups ${J(X)}$. {IV}.
\newblock {\em Topology}, 5:21--71, 1965.

\bibitem[And]{weakmaps}
L.~Anderson.
\newblock Representing weak maps of oriented matroids.
\newblock To appear in {\em European Journal of Combinatorics}.

\bibitem[And98]{homotopy}
L.~Anderson.
\newblock Homotopy groups of the combinatorial {G}rassmannian.
\newblock {\em Discrete Comput. Geom.}, 20:549--560, 1998.

\bibitem[And99a]{AB}
L.~Anderson.
\newblock Matroid bundles.
\newblock In {\em New Perspectives in Algebraic Combinatorics}, MSRI book
  series. Cambridge University Press, 1999.

\bibitem[And99b]{CDmanifolds}
L.~Anderson.
\newblock Topology of combinatorial differential manifolds.
\newblock {\em Topology}, 38(1):197--221, 1999.

\bibitem[Bab93]{Ba}
E.~Babson.
\newblock {\em A combinatorial flag space}.
\newblock PhD thesis, MIT, 1993.

\bibitem[BLS{\etalchar{+}}93]{BLSWZ}
A.~Bj{\"o}rner, M.~{Las Vergnas}, B.~Sturmfels, N.~White, and G.~M. Ziegler.
\newblock {\em Oriented matroids}, volume~46 of {\em Encyclopedia of
  Mathematics and its Applications}.
\newblock Cambridge University Press, 1993.

\bibitem[Dol63]{D}
A.~Dold.
\newblock Partitions of unity in the theory of fibrations.
\newblock {\em Ann. of Math. (2)}, 78:223--255, 1963.

\bibitem[DT56]{DT}
A.~Dold and R.~Thom.
\newblock Une g\'en\'eralisation de la notion d'espace fibr\'e. {A}pplication
  aux produits sym\'etriques infinis.
\newblock {\em C. R. Acad. Sci. Paris}, 242:1680--1682, 1956.

\bibitem[FL78]{FL}
J.~Folkman and J.~Lawrence.
\newblock Oriented matroids.
\newblock {\em J. Combin. Theory Ser. B}, 25:199--236, 1978.

\bibitem[GM92]{GM}
I.~M. Gelfand and R.~D. MacPherson.
\newblock A combinatorial formula for the {P}ontrjagin classes.
\newblock {\em Bull. Amer. Math. Soc. (N.S.)}, 26:304--309, 1992.

\bibitem[Hir75]{Hiro}
H.~Hironaka.
\newblock Triangulations of algebraic sets.
\newblock In {\em Algebraic geometry (Proc. Sympos. Pure Math., Vol. 29,
  Humboldt State Univ., Arcata, Calif., 1974)}, pages 165--185. Amer. Math.
  Soc., 1975.

\bibitem[Mac93]{Mac}
R.~D. MacPherson.
\newblock Combinatorial differential manifolds.
\newblock In {\em Topological methods in modern mathematics (Stony Brook, NY,
  1991)}, pages 203--221. Publish or Perish, 1993.

\bibitem[Mil59]{M}
J.~Milnor.
\newblock On spaces having the homotopy type of a {CW}-complex.
\newblock {\em Trans. Amer. Math. Soc.}, 90:272--280, 1959.

\bibitem[MRG93]{MR93}
N.~Mn{\"e}v and J.~Richter-Gebert.
\newblock Two constructions of oriented matroids with disconnected extension
  space.
\newblock {\em Disc. and Comp. Geometry}, 10(3):271--286, 1993.

\bibitem[MS74]{MS}
J.~W. Milnor and J.~D. Stasheff.
\newblock {\em Characteristic classes}.
\newblock Number~76 in Annals of Mathematics Studies. Princeton University
  Press, 1974.

\bibitem[MZ93]{MZ}
N.~Mn{\"e}v and G.~Ziegler.
\newblock Combinatorial models for the finite-dimensional {G}rassmannians.
\newblock {\em Discrete Comput. Geom.}, 10(3):241--250, 1993.

\bibitem[Qui73]{Q73}
D.~Quillen.
\newblock Higher algebraic ${K}$-theory, {I}: Higher ${K}$-theories.
\newblock In {\em Proc. Conf., Battelle Memorial Inst., Seattle, Wash., 1972},
  number 341 in Lecture Notes in Mathematics, pages 85--147. Springer-Verlag,
  1973.

\bibitem[Sta63]{S}
J.~Stasheff.
\newblock A classification theorem for fibre spaces.
\newblock {\em Topology}, 2:239--246, 1963.

\bibitem[SZ93]{SZ}
B.~Sturmfels and G.~Ziegler.
\newblock Extension spaces of oriented matroids.
\newblock {\em Discrete Comput. Geom.}, 10(1):23--45, 1993.

\bibitem[Whi78]{Whitehead}
G.~W. Whitehead.
\newblock {\em Elements of Homotopy Theory}.
\newblock Springer-Verlag, 1978.

\end{thebibliography}
